\documentclass[11pt]{article}
\usepackage[latin1]{inputenc}
\usepackage{amssymb}
\usepackage{latexsym}
\usepackage{color}
\usepackage{amsmath}
\usepackage{amsfonts}
\usepackage{amsthm}
\usepackage{varioref}
\usepackage{enumerate,pdfsync}
\usepackage{cancel}

\def\be{\begin{eqnarray}}
\def\ee{\end{eqnarray}}
\def\beq{\begin{equation}}
\def\eeq{\end{equation}}
\def\b*{\begin{eqnarray*}}
\def\e*{\end{eqnarray*}}
\def\bi{\begin{itemize}}
\def\ei{\end{itemize}}

\def \1{{\bf 1}}

\def\eps{\varepsilon}

\def\={\;=\;}

\def\x{\times}

\def\Pro#1{\mathbb{P}\left[{#1}\right]}

\def\essinf{{\rm ess}\!\inf\limits}
\def\esssup{{\rm ess}\!\sup\limits}

\def \proof{{\noindent \bf Proof. }}
\def \ep{\hbox{ }\hfill$\Box$}
 \def\reff#1{{\rm(\ref{#1})}}

\def\Pas{\mathbb{P}-\mbox{a.s.}}

 \def\vs#1{\vspace{#1mm}}
 
 \def\no{\noindent}

\def \F{\mathbb{F}}

\def \H{\mathbb{H}}

\def \N{\mathbb{N}}
\def \P{\mathbb{P}}
\def \Q{\mathbb{Q}}
\def \R{\mathbb{R}}

\def\Ec{{\cal E}}
\def\Fc{{\cal F}}
\def\Gc{{\cal G}}

\def\Jc{{\cal J}}

\def\Kc{{\cal K}}

\def\Tc{{\cal T}}

\def\Xc{{\cal X}}
\def\Yc{{\cal Y}}
\def\Zc{{\cal Z}}

\def\Lb{{\mathbf L}}
\def \Hb{{\mathbf H}}

\newtheorem{Theorem}{Theorem}[part]

\newtheorem{Proposition}{Proposition}[part]

\newtheorem{Lemma}{Lemma}[part]
\newtheorem{Corollary}{Corollary}[part]
\newtheorem{Remark}{Remark}[part]

\makeatletter \@addtoreset{equation}{section}

\@addtoreset{Definition}{section}

\@addtoreset{Theorem}{section}

\@addtoreset{Proposition}{section}

\@addtoreset{Assumption}{section}

\@addtoreset{Corollary}{section}

\@addtoreset{Lemma}{section}

\@addtoreset{Remark}{section}

\@addtoreset{Example}{section}

\def\Indi#1{{\rm \bf 1}_{#1}}
\def\dom#1{{\rm dom}(#1)}
\def\Leb{{\rm Leb}}

\theoremstyle{plain}
\newtheorem{theorem}{Theorem}[section]

\newtheorem{definition}[theorem]{Definition}
\theoremstyle{definition}

\def\Ab{{\bf A}}
\def\Sb{{\bf S}}
\def\Kb{{\bf K}}

\def\Lambdab{{\bf   \Lambda}}

\def\vc{\vartheta}

\def\Pae{\P-\mbox{a.e.~}}

\def\real{\mathbb{R}}
\def\BSDE{{\rm BSDE}}
\def\Eg#1#2{\Ec^{g}_{#1} [#2]}
\def\H#1{{\rm\bf (H$_{#1}$)}}
\def\Hgconv{{\rm\bf (H$_{conv}$)}}
\def\HPsi{\H{\Psi}}
\def\Hg{\H{g}}

\title{BSDEs with weak terminal condition}
\author{Bruno Bouchard\footnote{and CREST, bouchard@ceremade.dauphine.fr} \and Romuald Elie\footnote{elie@ceremade.dauphine.fr} \and Anthony R\'eveillac\footnote{anthony.reveillac@ceremade.dauphine.fr} \and \\ Universit\'e Paris-Dauphine\\ CEREMADE UMR CNRS 7534 \\Place du Mar\'echal De Lattre De Tassigny \\ 75775 Paris cedex 16 France
}

\begin{document}
\maketitle 
\begin{abstract} We introduce a new class of   Backward Stochastic Differential Equations in which the $T$-terminal value $Y_{T}$ of the solution $(Y,Z)$ is not fixed as a random variable, but 
only satisfies a weak constraint of the form $E[\Psi(Y_{T})]\ge m$, for some (possibly random) non-decreasing map $\Psi$ and some threshold $m$. We name them \textit{BSDEs with weak terminal condition} and obtain a  representation of the minimal time $t$-values $Y_{t}$ such that $(Y,Z)$ is a supersolution of the BSDE with weak terminal condition.
It provides  a non-Markovian BSDE formulation of the PDE characterization  obtained   for Markovian stochastic target problems under controlled loss  in Bouchard, Elie and Touzi \cite{BoElTo09}. We then study the main properties of this minimal value. In particular, we analyze its continuity and convexity with respect to the $m$-parameter appearing in the weak terminal condition, and show how it can be related to a dual optimal control problem in Meyer form. These last properties generalize to a non Markovian framework previous results on quantile hedging and hedging under loss constraints obtained in   F\"{o}llmer and Leukert \cite{FoLe99,FoLe00}, and in  Bouchard, Elie and Touzi \cite{BoElTo09}. 
\end{abstract}

\vspace{5mm}


\vspace{5mm}

\noindent {\textbf{MSC Classification (2000)}: Primary: 60H10; 93E20; Secondary: 49L20; 91G80 
\newpage

\section{Introduction}

Solving a backward stochastic differential equation (hereafter BSDE), with terminal data $\xi \in \Lb_2(\Fc_{T})$ and driver $g$, consists in finding a pair of predictable processes $(Y,Z)$, with certain integrability properties, such that the dynamics of $Y$ satisfies 
$
dY_{t}=-g(t,Y_{t},Z_{t})dt + Z_{t} dW_{t}
$
 and 
 $
 Y_{T}=\xi
 $ (where $W$ denotes a standard Brownian motion).
 It can  be rephrased in: find an  initial data $Y_{0}$ and a control process $Z$ such that the solution  $Y^{Z}$ of the controlled stochastic differential equation  
 \be\label{eq: intro def BSDE as controlled sde}
 Y^{Z}_{t}=Y_{0}-\int_{0}^{t}  g(s,Y^{Z}_{s},Z_{s})ds + \int_{0}^{t} Z_{s} dW_{s} \;, \qquad 0\le t \le T\;,
 \ee
 satisfies $Y^{Z}_{T}= \xi$.   In cases where the previous problem does not admit a solution, a weaker formulation is to find an initial data $Y_{0}$ and a control $Z$  such that  
\be\label{eq: intro contrainte super hedge}
Y^{Z}_{T}\ge  \xi\;\;\Pas
\ee 
In most applications, one is interested in the minimal initial condition $Y_{0}$ and in the associated control $Z$. This is for instance the case in the financial literature  in which $Y_{0}$ represents the cost of the cheapest super-replication strategy for the contingent claim $\xi$, and $Z$ provides the associated hedging strategy, see e.g. \cite{KaPeQu97}.
 
Motivated by situations where this minimal  value $Y_{0}$ is too large for practical applications, it was suggested to relax the strong constraint   \reff{eq: intro contrainte super hedge} into a weaker one of the form \be\label{eq: intro contrainte moment}
E\left[\ell(Y^{Z}_{T}-\xi)\right]\ge m \;,
\ee 
where $m$ is a given threshold and $\ell$ is a non-decreasing map. For $\ell(x)=\1_{\{x\ge 0\}}$, this corresponds to matching the criteria $Y^{Z}_{T}\ge  \xi$ at least with probability $m$. In financial terms, this is  the so-called quantile hedging problem, see  \cite{FoLe99}\footnote{ {In fact, their original formulation also imposes a  budget constraint constraint $Y^{Z}_{T}\ge 0$ $\Pas$, which can be taken into account by imposing a criteria of the form \eqref{eq:psi1} with    $\Psi(Y^{Z}_{T}):=\1_{\{Y^Z_T-\xi\ge 0\}} - \infty \1_{\{Y_T^Z<0\}}$.}} . More generally, $\ell$ is viewed as a loss function, one typical example being $\ell(x):=-(x^{-})^{q}$ with $q\ge 1$, see   \cite{FoLe00} for general non-Markovian but linear dynamics. Such problems were coined ``stochastic target problems with controlled loss'' by  \cite{BoElTo09} who consider  a non-linear  Markovian formulation in a Brownian diffusion setting, see also   \cite{Mo10} for the jump diffusions setting.

The aim of this paper is to study the non-linear non-Markovian setting in which the terminal constraint is of the form 
\be
\label{eq:psi1}
 &E\left[\Psi(Y^{Z}_{T})\right]\ge m.& \label{eq: intro contrainte moment psi}
\ee
In the above, $m\in \R$ and $\Psi$ is a (possibly random) non-decreasing real-valued map.  Our problem can then be written as 
\be\label{eq: intro def bsde weak terminal cond}
&\mbox{ Find the minimal $Y_{0}$ such that \reff{eq: intro def BSDE as controlled sde} and \reff{eq: intro contrainte moment psi}   hold  for some $Z$.} &
\ee 
This leads to the introduction of a new class of BSDEs which we call BSDEs with weak terminal condition. More precisely,
we refer to this problem by saying that we want to solve the BSDE with driver $g$ and weak terminal condition $(\Psi,m)$ to insist on the fact that the terminal condition $Y^{Z}_{T}$ is not fixed as a random variable, but only has to satisfy the weak constraint \reff{eq: intro contrainte moment psi}.
  
The first step in our analysis lies in a reformulation based on the martingale representation theorem, as suggested in \cite{BoElTo09}. More precisely, if $Y_{0}$ and $Z$ are such that \reff{eq: intro contrainte moment psi} holds, then the martingale representation Theorem implies that we can find an element $\alpha$ in the set $\Ab_0$,  of predictable   square integrable processes, such that 
\b*
 \Psi(Y^{Z}_{T})\ge M^{\alpha}_{T}:=m+\int_{0}^{T} \alpha_{s} dW_{s}.
 \e*
 On the other hand,   since $\Psi$ is non-decreasing, one can introduce its  {left-continuous} inverse $\Phi$ and note that the solution $(Y^{\alpha},Z^{\alpha})$ of the BSDE
 \begin{equation}\label{eq: intro BSDE alpha}
 Y^{\alpha}_{t}=\Phi(M^{\alpha}_{T})+\int_{t}^{T} g(s,Y^{\alpha}_{s},Z^{\alpha}_{s})ds -\int_{t}^{T} Z^{\alpha}_{s} dW_{s} ,\; 0\le t \le T\;,
 \end{equation}
 actually solves \reff{eq: intro def BSDE as controlled sde} and \reff{eq: intro contrainte moment psi}. We indeed show that the solution of \reff{eq: intro def bsde weak terminal cond} is given by
 \be\label{eq: inf alpha Y0 alpha}
 \inf\{Y_{0}^{\alpha},\;\alpha \in \Ab_0\}.
 \ee   
 This leads to study its dynamical counterpart  
 \begin{equation}\label{introeq1}
 \Yc^{\alpha}_{\tau}:=\essinf\{Y^{\alpha'}_{\tau}, \;\alpha'\in \Ab_0 \mbox{ s.t. } \;\alpha'=\alpha \mbox{ on } [\![0,\tau]\!]\}\;, \quad 0\le \tau \le T\;.
 \end{equation}
 We verify that the family $\{ \Yc^{\alpha},\;\alpha\in \Ab_0\}$ satisfies a dynamic programming principle which can be seen as a counterpart of the geometric dynamic programming principle of   \cite{SoTo02} used in \cite{BoElTo09}.  In particular, this implies that $\{ \Yc^{\alpha},\;\alpha\in \Ab_0\}$  is a $g$-submartingale family to which we can apply the non-linear Doob-Meyer decomposition of  \cite{Pe99}. This provides a representation of the family $\{ \Yc^{\alpha},\;\alpha\in \Ab_0\}$ in terms of minimal supersolutions to a family of BSDEs with driver $g$ and (strong) terminal conditions $\{ \Phi(M^{\alpha}_T),\;\alpha\in \Ab_0\}$. This representation allows in particular to characterize the family $\{ \Yc^{\alpha},\;\alpha\in \Ab_0\}$ uniquely. Under additional convexity assumptions on the coefficients $g$ and $\Phi$, we observe that the essential infimum in \reff{introeq1} is attained. Hence, there exists an optimal $\hat\alpha\in\Ab_0$ such that solving the BSDE with weak terminal condition $(\Psi,m)$ boils down to solving the BSDE with dynamics \reff{eq: intro BSDE alpha} and strong terminal condition $\Phi(M^{\hat\alpha}_{T})$. In a Markovian framework, our approach provides in particular a BSDE formulation for the PDEs derived in \cite{BoElTo09}.
\\
 We then study in details important properties of this family and focus in particular on the regularity of $\Yc^{\alpha}$ with respect to the threshold parameter $m$.
  We exhibit, under weak conditions, a stability property of the solution with respect to the variations of the parameter $m$. We also observe  that $\Yc^{\alpha}$  is convex with respect to the threshold parameter. This observation allows us in particular to conclude that $\Phi$ (whenever it is deterministic) can be replaced by its 
more regular convex envelope in order to compute $\Yc^{\alpha}$ on $[0,T)$. This was already observed in the restrictive Markovian setting of  \cite{BoElTo09}, in which it is proved by using PDE technics. We provide here a pure probabilistic argument. Similarly, it was also observed in \cite{FoLe99}, \cite{FoLe00} and  \cite{BoElTo09} that \reff{eq: intro def bsde weak terminal cond} admits a dual linear problem when $g$ is linear. We extend this result via probabilistic arguments to the semi-linear setting, for which the dual formulation takes the form of a stochastic control problem in Meyer form.\\ 

The rest of the paper is organized as follows. In Section \ref{SEC: BSDE weak terminal condition}, we provide a precise formulation for \reff{eq: intro def bsde weak terminal cond} and relate this problem to a $g$-submartingale family satisfying a dynamic programming principle. Attainability of the optimal control $\hat \alpha\in\Ab_0$ is also discussed. Section \ref{SEC: Properties} collects the continuity and convexity properties as well as the dual formulation of the problem. 
Finally, Section \ref{SEC: proof thm BSDE characterization bar Yc} contains the proof of the BSDE representation for $\{ \Yc^{\alpha},\;\alpha\in \Ab_0\}$.\\

We close this introduction with a series of notations that will be used all over this paper.  Let $d\geq 1$ and $T>0$ be fixed. We denote by $W:=(W_t)_{t\in [0,T]}$ a $d$-dimensional Brownian motion defined on a probability space $(\Omega,\mathcal{F},\P)$ with $\P$-augmented natural filtration $\F=(\mathcal{F}_t)_{t\in [0,T]}$. The components of $W$ will be denoted by $W=(W^1,\cdots,W^d)$ and $E$ will stand for the expectation with respect to $\P$. For simplicity, we assume that $\mathcal{F}=\mathcal{F}_T$.  Throughout the paper we will make use of the following spaces. 
\begin{itemize}
\item[-] $\Lb_{p}(U,\Gc)$ denotes the set of $p$-integrable $\Gc$-measurable random variables with values in $U$, $p\ge 0$, $U$ a Borel set of $\R^{n}$ for some $n\ge 1$ and $\Gc\subset \Fc$. When $U$ and $\Gc$ can be clearly identified by the context, we omit them. This will be in particular the case when $\Gc=\Fc$.
\item[-] $\Tc$ denotes the set of $\F$-stopping times in $[0,T]$. For $\tau_{1}\in \Tc$, $\Tc_{\tau_{1}}$ is the set of stopping times $\tau_{2}$ in $\Tc$ such that $\tau_{2}\ge \tau_{1}$ $\Pas$ The notation $E_{\tau}[\cdot]$ stands for the conditional expectation given $\Fc_{\tau}$, $\tau \in \Tc$.
\item[-] $\Sb_{2}$ denotes the set of $\real$-valued, \textit{c\`adl\`ag}\footnote{right-continuous with left limits} and $\F$-adapted stochastic processes $X=(X_t)_{t\in [0,T]}$ such that $\|X\|_{\Sb_{2}}:=E[\sup_{t \in [0,T]}|X_t|^2]^{1/2}<\infty$.
\item[-]   $\Hb_{2}$ denotes the set of $\real^{n}$-valued, $\F$-predictable stochastic processes $X=(X_t)_{t\in [0,T]}$ such that $\|X\|_{\Hb_{2}}:=E\left[\int_0^T |X_t|^2 dt\right]^{1/2}<\infty$. In the following, the dimension $n$ will be given by the context. 
\item[-] $\Kb_{2}$ denotes the set of non-decreasing $\real$-valued and $\F$-adapted stochastic processes $X=(X_t)_{t\in [0,T]}$ such that $\|X\|_{\Sb_{2}}<\infty$.
\end{itemize} 
 
 Inequalities between random variables are understood in the $\Pas$-sense.
\section{BSDE with weak terminal condition}\label{SEC: BSDE weak terminal condition}

\subsection{Definitions and problem reformulation}

We first define the main object of this paper. 

\begin{definition}[Solution to a  BSDE with weak terminal condition] 
Given a measurable map $\Psi: \R\x \Omega \mapsto U$, with $U\subset \R\cup \{-\infty\}$, $\tau \in \Tc$ and $\mu\in \Lb_0(U,\Fc_{\tau})$, we say that $(Y,Z)\in \Sb_{2}\x \Hb_{2}$ is a supersolution of  the BSDE with generator $g:\Omega \times [0,T]\times \real \times \real^d \to \real $ and  weak terminal condition $(\Psi,\mu,\tau)$, in short $\BSDE(g,\Psi,\mu,\tau)$, if for any  {$0\leq s\leq t\leq T$},
\begin{align}
\label{eq: def BSDEw(psi,p)} &Y_s \ge Y_{t} +\int_{s}^{t} g(r,Y_{r},Z_{r}) dr- \int_{s}^{t} Z_{r} dW_{r}, \quad \mbox{ and }\\
\label{eq: def BSDEw(psi,p)term} &E_{\tau}\left[\Psi(Y_{T})\right]\ge \mu.
\end{align}

\end{definition}

\noindent
Before discussing the well-posedness of Equation \reff{eq: def BSDEw(psi,p)}-\reff{eq: def BSDEw(psi,p)term}, let us emphasize that the difference with classical BSDEs lies in the fact that we do not prescribe a terminal condition to $Y$ in the classical $\Pas$-sense but only impose a weak condition in expectation form (which justifies the terminology of \textit{BSDE with weak terminal condition}). Even if we were asking for equalities in  \reff{eq: def BSDEw(psi,p)}-\reff{eq: def BSDEw(psi,p)term}, this would obviously be too weak to expect uniqueness, as any random variable $\xi$ satisfying $E_{\tau}\left[\Psi(\xi)\right]= \mu$ could serve as a terminal condition.  

However, when $\Psi$ is non-decreasing, the set 
\begin{equation}\label{eq: def Gamma}
\Gamma(\tau,\mu):=\{Y_{\tau}: (Y,Z)\in \Sb_{2}\x \Hb_{2} \mbox{ is a supersolution of } \BSDE(g,\Psi,\mu,\tau)\} \;,
\end{equation}
defined for any $\tau\in\Tc$ and $\mu\in\Lb_0(U,\Fc_{\tau})$, can be characterized by its lower-bound, whenever it is achieved. \\

Throughout the paper, we shall restrict to the  case where $g$ is Lipschitz continuous with linear growth, $\Psi^{+}$ is bounded, and  the domain of $\Psi$ is bounded from below,  in order to avoid un-necessary technicalities. 
\vs1

\textbf{Standing Assumption \HPsi:} For $\Pae$ $\omega \in \Omega$, the map $y\in \R \mapsto \Psi( \omega,y)$ is non-decreasing,  {right-continuous}, valued in $[0,1]\cup\{-\infty\}$, and  its left-continuous inverse 
$\Phi(\omega,\cdot)$ satisfies  
$\Phi:\Omega \times [0,1] \mapsto [0,1]$  is measurable. 
\vs1

  By left-continuous inverse we mean  the left-continuous map defined for $\omega$ fixed by
$$\Phi(\omega,x):=\inf\{y\in \real,\; \Psi(\omega,y)\geq x\},$$ which satisfies
 \begin{equation}
\label{eq:relationsphipsi}
\Phi \circ \Psi \leq \textrm{Id}\le \Psi \circ \Phi. 
\end{equation} 
 The left-hand side follows from the definition of $\Phi$, the right-hand side  holds by right-continuity of $\Psi$.  
Note that the above assumption implies   $\Psi(\omega,\cdot)=-\infty$ on $(-\infty,0)$ and $\Psi(\omega,\cdot)=1$ on $[1,\infty)$.  In particular,  the constraint in expectation \reff{eq: def BSDEw(psi,p)term} implies $Y_{T}\ge 0$ $\Pas$ Obviously the set $[0,1]$ is chosen for ease of notations and can be replaced by any closed interval.

\vs1
\textbf{Standing Assumption \Hg}  $g$ is a measurable map from $\Omega \times [0,T] \times \real \times \real^d$ to $\real$ and  $g(\cdot,y,z)$ is $\F$-predictable, for each $(y,z)\in \R\x \R^{d}$. There exists a constant $K_{g}>0$ and a random variable  $\chi_{g}\in   \Lb_{2}(\R_{+})$,  such that  
\b* 
&|g(t,0,0)|\le \chi_{g} \; \;\Pas&\\
&|g(t,y_1,z_1)-g(t,y_2,z_2)| \leq  K_{g} (|y_1-y_2| +|z_1-z_2|)\; \;\Pas&
\\
&\forall (t,y_i,z_i) \in [0,T]\times \real \times \real^d, \; i=1,2.&
\e*   

Let $\Ab_{\tau,\mu}$ denote the set elements $\alpha \in \Hb_{2}$ such that 
\be\label{eq: def P alpha}
M^{(\tau,\mu),\alpha}  := \mu+ \int_{\tau}^{\tau \vee \cdot} \alpha_{s} dW_{s} \;\mbox{ takes values in } [0,1].  
\ee
Then, \eqref{eq: def BSDEw(psi,p)term} is equivalent to   $\Psi(Y_T) \ge M_T^{(\tau,\mu),\alpha}$ for some $\alpha \in \Ab_{\tau,\mu}$. In view of \eqref{eq:relationsphipsi}, this is equivalent to $Y_{T}\ge \Phi(M^{(\tau,\mu),\alpha}_T)$ for some $\alpha \in \Ab_{\tau,\mu}$. This implies that supersolutions  of  BSDE$(g,\Psi,\mu,\tau)$ can be characterized in terms of $g$-expectations whose definition is recalled below. 

\begin{definition}[g-expectation] \label{def: g expectation} Given $\tau_{2}\in \Tc$ and  $\xi \in \Lb_2(\R,\Fc_{\tau_{2}})$, let $(Y,Z) \in \Sb_{2}\x \Hb_{2}$ denote the solution of 
\b*
Y=\xi+ \int_{\cdot\wedge \tau_{2}}^{\tau_{2}} g(s,Y_{s},Z_{s})  {ds}- \int_{\cdot\wedge \tau_{2}}^{\tau_{2}} Z_{s} dW_{s}.
\e* 
Then, we define the (conditional) $g$-expectation of $\xi$ at the stopping time $\tau_{1}\le \tau_{2}$ as $\Eg{\tau_{1},\tau_{2}}{\xi}:=Y_{\tau_{1}}$. When $\tau_{2}\equiv T$, we only write 
 $\Eg{\tau_{1}}{\xi}$, and say that $(Y,Z)$ solves $\BSDE(g,\xi)$. 
\end{definition}
 
Note that existence and uniqueness hold under Assumption \Hg. 
In the following, we shall adopt the terminology of  Peng \cite{Pe04} and call $g$-martingale  (resp. $g$-submartingale) a  process $Y$ such that $\Eg{t,s}{Y_{s}}=Y_{t}$ (resp. $\Eg{t,s}{Y_{s}}\ge  Y_{t}$), for all $t\le s\le T$.

\begin{Proposition}\label{prop: lien BSDE(psi,p) et BSDE(alpha)}  Fix $\tau \in \Tc, \;\mu\in \Lb_0([0,1],\Fc_{\tau})$. Then,  
$(Y,Z) \in \Sb_{2} \times \Hb_{2}$ is a supersolution of  $\BSDE(g,\Psi,\mu,\tau)$ if and only if $(Y,Z)$ satisfies \reff{eq: def BSDEw(psi,p)} and there exists $\alpha\in \Ab_{\tau,\mu}$ such that $Y_{t}\ge \Eg{t}{\Phi(M^{(\tau,\mu),\alpha}_T)}$ for $t\in [0,T]$ $\Pas$
\end{Proposition}

\proof 
Let $(Y,Z)$ be a super solution  of BSDE$(g,\Psi,\mu,\tau)$. Then there exists some element $\rho$ in $\Lb_0([0,1],\Fc_{\tau})$ with $\rho\geq \mu$, $\Pas$ and $\tilde{\alpha}$ in $\Ab_{\tau,\rho}$ such that $\Psi(Y_T)=M_T^{(\tau,\rho),\tilde{\alpha}}$. Set $\theta^{\tilde{\alpha}}:=\inf\{s\geq \tau, \; M_s^{(\tau,\mu),\tilde{\alpha}}=0\}$. It is clear that $\theta^{\tilde{\alpha}}$ belongs to $\Tc$ and that $\alpha:=\tilde{\alpha} \textbf{1}_{[ 0,\theta^{\tilde{\alpha}})}$ belongs to {$\Ab_{\tau,\mu}$} and satisfies $M_T^{(\tau,\rho),\tilde{\alpha}} \geq M_T^{(\tau,\mu),\alpha}$, $\Pas$, since $M_T^{(\tau,\rho),\tilde{\alpha}} \geq 0$ by definition of $\Ab_{\tau,\rho}$. The monotonicity of $\Phi$ and   \eqref{eq:relationsphipsi} imply that 
$$ Y_T \geq (\Phi \circ \Psi)(Y_T) \geq \Phi(M_T^{(\tau,\mu),\alpha}). $$
By comparison for Lipschitz BSDEs, we obtain $Y_{t}\ge \Eg{t}{\Phi(M^{(\tau,\mu),\alpha}_T)}$ for $t\in [0,T]$. Conversly, let $\alpha\in \Ab_{\tau,\mu}$ be such that $Y_{t}\ge \Eg{t}{\Phi(M^{(\tau,\mu),\alpha}_T)}$ for $t\in [0,T]$ and assume that $(Y,Z)$ satisfies \reff{eq: def BSDEw(psi,p)}. Then,   \reff{eq:relationsphipsi} implies
$$ \Psi(Y_T) \geq (\Psi\circ \Phi)(M^{(\tau,\mu),\alpha}_T) \ge M^{(\tau,\mu),\alpha}_T.$$
Taking the conditional expectation on both sides leads to \eqref{eq: def BSDEw(psi,p)term}. 
\ep 
\vs2
 
  In view of Proposition \ref{prop: lien BSDE(psi,p) et BSDE(alpha)}, the lower bound of $\Gamma(\tau,\mu)$ (which we recall, has been defined in \eqref{eq: def Gamma}) can be expressed in terms of 
 \be\label{eq : def bar Y}
~~~\Yc_{\tau}(\mu):= \essinf_{\alpha\in \Ab_{\tau,\mu}} \Ec^g_{\tau}\left[\Phi(M^{(\tau,\mu),\alpha}_T)\right]\;,\;\tau \in \Tc, \;\mu\in \Lb_{0}([0,1],\Fc_{\tau}).
\ee 
 
This is the statement of the next corollary. 

\begin{Corollary}
$\essinf\Gamma(\tau,\mu)=\Yc_{\tau}(\mu)$,   $\forall$ $\tau \in \Tc$, $\mu\in \Lb_0([0,1],\Fc_{\tau})$.
\end{Corollary}

\proof  The fact that $Y_{\tau} \in\Gamma(\tau,\mu)$ implies $Y_{\tau} \ge \Yc_{\tau}(\mu)$ follows from Proposition \ref{prop: lien BSDE(psi,p) et BSDE(alpha)}.
On the other hand, the same proposition implies that each $\Eg{\tau}{\Phi(M^{(\tau,\mu),\alpha}_T)}$ with $\alpha\in \Ab_{\tau,\mu}$ belongs to $\Gamma(\tau,\mu)$.\ep

 \begin{Remark}\label{rem: borne unif Yc}{\rm For later use, note that the assumptions \Hg~and \HPsi~ensure that we can find $\eta \in \Sb_{2}$ such that $|\Ec^{g}_{t}[\Phi(M)]|\vee |\Yc_{t}(\mu)|\le \eta_{t}$,  for all $t\le T$ and $\mu\in \Lb_{0}([0,1],\Fc_{t})$, $M\in \Lb_{0}([0,1])$. See (i) of Proposition \ref{prop: majo bsde appendix} in the Appendix.}
 \end{Remark}
 
 \begin{Remark}\label{rem: conca sur event set}{\rm Note that $\Yc_{\tau}(\mu)= \Yc_{\tau}(\mu_{1})\1_{A}+\Yc_{\tau}(\mu_{2})\1_{A^{c}}$ whenever $\mu:=\mu_{1}\1_{A}+\mu_{2}\1_{A^{c}}$ for    $A\in \Fc_{\tau}$, $\mu_{1},\mu_{2}\in \Lb_0([0,1],\Fc_{\tau})$, and $\tau \in \Tc$.  
Indeed, $\alpha:=\1_{[\tau,T]}(\alpha_{1}\1_{A}+\alpha_{2}\1_{A^{c}})\in  \Ab_{\tau,\mu}$  for all $\alpha_{i}\in  \Ab_{\tau,\mu_{i}}$ with $i=1,2$. Since $\Ec^g_{\tau}\left[\Phi(M^{(\tau,\mu),\alpha}_T)\right]=\Ec^g_{\tau}\left[\Phi(M^{(\tau,\mu_{1}),\alpha_{1}}_T)\right]\1_{A}+ \Ec^g_{\tau}\left[\Phi(M^{(\tau,\mu_{2}),\alpha_{2}}_T)\right]\1_{A^{c}}$, this implies  $\Yc_{\tau}(\mu)\le \Yc_{\tau}(\mu_{1})\1_{A}+\Yc_{\tau}(\mu_{2})\1_{A^{c}}$. The converse inequality follows from the previous identity applied with $\alpha_{1}:=\alpha\1_{A}$ and $\alpha_{2}:=\alpha\1_{A^{c}}$ for any $\alpha\in \Ab_{\tau,\mu}$ so that $\alpha_{i}\in  \Ab_{\tau,\mu_{i}}$ for $i=1,2$.
 }
 \end{Remark} 
 \begin{Remark} {\rm  Before going on with the study of the set $\Gamma$, let us notice that a similar analysis  can be carried out for weak constraints of the form $\Ec^{h}_{\tau}\left[\Psi(Y_{T})\right]\ge \mu$ in place of $E_{\tau}\left[\Psi(Y_{T})\right]\ge \mu$ in \reff{eq: def BSDEw(psi,p)term}, with $\Ec^h$ defined as the $h$-expectation associated to some random map $h$ satisfying similar conditions as $g$. In finance, the latter condition interprets as a risk-measure constraints, see e.g. \cite{Pe04}, while our condition is more related to expected loss constraints, see \cite{FoLe00}. Again, we try to avoid un-necessary additional technicalities and stick to the case $h\equiv0$.  
   }
 \end{Remark}
 
\subsection{BSDE characterization of the minimal initial condition}\label{SEC: main characterization}
 
The main result of this section is a BSDE characterization for the lower bound of the set  $\Gamma(\tau,\mu)$ of time-$\tau$ initial conditions of supersolutions of $\BSDE(g,\Psi,\mu,\tau)$. In particular, this extends to a non Markovian framework the PDE characterization of \cite{BoElTo09}.
 
\vs2

For ease of notations, we now fix $m_{o}\in [0,1]$ and set
$$
\left\lbrace
\begin{array}{l}
M^{\alpha}_t:=M_t^{(0,m_{o}),\alpha} \mbox{ , } 
\Ab^{\alpha}_{\tau}:=\{\alpha' \in \Ab_{\tau,M^\alpha_\tau}:~\alpha'=\alpha\;dt\times d\P \mbox{ on } [\![0,\tau]\!]\},\\
\Ab_{0}:=\Ab_{0,m_{o}} \mbox{ and } \Yc^{\alpha}_t:=\Yc_t (M^{\alpha}_t) \mbox{ for } \alpha \in \Ab_{0}, \; t \in [0,T],
\end{array}
\right.
$$
where we recall that $M^{(0,m_{o}),\alpha}$ and $\Ab_{0,m_{o}}$ are given in \eqref{eq: def P alpha}.
 
\begin{Theorem}\label{thm : BSDE characterization bar Yc} ~ For any $\alpha \in \Ab_{0}$, $\Yc^{\alpha}$ is a  $g$-submartingale,  {it is l\`adl\`ag \footnote{left and right-limited according to the french celebrated acronym} on countable sets}, and the following dynamic programming principle holds:

\begin{itemize}

\item[(i)]\label{item: dpp} $\Yc^{\alpha}_{\tau_{1}}=\essinf_{\bar \alpha \in \Ab^{\alpha}_{\tau_{1}}}\Eg{\tau_{1},\tau_{2}}{\Yc^{\bar \alpha}_{\tau_{2}}}$,  for each $\tau_{1}\in \Tc$, $\tau_{2}\in \Tc_{\tau_{1}}$. 

\end{itemize}
Under the additional assumption that 
\begin{equation}
\label{eq:hypcontiPhi}
m\in [0,1]\mapsto \Phi(\omega,m) \mbox{ is continuous for } \P\mbox{-a.e. } \omega \in \Omega, 
\end{equation}
the following holds:
\begin{itemize}

\item[(ii)]\label{item : indistingable proc cadlag}    $\Yc^{\alpha}$ is indistinguishable from a c\`adl\`ag $g$-submartingale, for each $\alpha \in \Ab_{0}$.

\item[(iii)] There exists a  family  $(  \Zc^{\alpha},   \Kc^{\alpha})_{\alpha \in \Ab_{0}}\subset  \Hb_{2}\x \Kb_{2}$  satisfying
 \be\label{eq : BSDE de la cond minimale - borne uniforme}
 \sup_{\alpha\in\Ab_{0}} \left\| (\Yc^\alpha, \Zc^\alpha, \Kc^\alpha)\right\|_{\Sb_{2}\times\Hb_{2}\times\Kb_{2}}
 &<& 
 \infty \;,
 \ee
 and such that, for all $\alpha \in \Ab_{0}$,  we have
\begin{equation}
\Yc^{\alpha}= \Phi(M^{\alpha}_{T})+\int_{\cdot}^{T}g(s,\Yc^{\alpha}_{s},\Zc^{\alpha}_{s})ds-\int_{\cdot}^{T} \Zc^{\alpha}_{s} dW_{s}+ \Kc^{\alpha} - \Kc^{\alpha}_T, \label{eq : BSDE de la cond minimale - BSDE} \\
\end{equation}
\begin{equation}
\Kc^{\alpha}_{\tau_{1}}= \essinf_{\bar\alpha\in \Ab_{\tau_{1}}^{\alpha} } E\left[\Kc^{\bar\alpha}_{\tau_{2}}|\Fc_{\tau_{1}}\right] \; ,\;\;\forall \;  \tau_{1}\in \Tc,\;\tau_{2}\in \Tc_{\tau_{1}},  \label{eq : BSDE de la cond minimale - cond de minimalite}
\end{equation} 
and 
\begin{equation}\label{eq : BSDE de la cond minimale - indep alpha futur}
(\Yc^{\alpha},\Zc^{\alpha},\Kc^{\alpha})\1_{[\![0,\tau]\!]}=(\Yc^{\bar\alpha},\Zc^{\bar\alpha},\Kc^{\bar\alpha})\1_{[\![0,\tau]\!]}, \;\;\forall \;    \tau\in \Tc,\; \bar\alpha \in \Ab^{\alpha}_{\tau}.
\end{equation}

\item[(iv)]  $(\Yc^{\alpha},\Zc^{\alpha}, \Kc^{\alpha})_{\alpha \in \Ab_{0}}$ is the unique family of $\Sb_{2}\x\Hb_{2}\x \Kb_{2}$ satisfying \reff{eq : BSDE de la cond minimale - borne uniforme}-\reff{eq : BSDE de la cond minimale - BSDE}-\reff{eq : BSDE de la cond minimale - cond de minimalite}-\reff{eq : BSDE de la cond minimale - indep alpha futur} for all $\alpha \in \Ab_{0}$.
\end{itemize}
\end{Theorem}
   
The proof of this theorem  is reported in Section \ref{SEC: proof thm BSDE characterization bar Yc}.

\begin{Remark} \label{rem: after theorem charact bsde}{\rm 

(i) The precise continuity assumption needed in the proof is :  $\Phi(M^{\alpha_{n}}_{T})$ converges in $\Lb_2$ to $\Phi(M^{\alpha}_{T})$ whenever $\|M^{\alpha_{n}}_{T}-M^{\alpha}_{T}\|_{\Lb_2}$ tends to $0$, for any sequence $(\alpha_{n})_{n}\subset \Ab_{0}$. However, this condition implies that $\Phi$ is continuous, as soon as random variables with non-absolutely continuous law with respect to the Lebesgue measure might be considered (which is the case here). 

(ii) We shall see in Proposition \ref{prop: regularisation au bord} below that $\Phi$ can  be replaced by its $m$-convex envelope, under mild assumptions. In this case, the   continuity assumption of the second part of Theorem \ref{thm : BSDE characterization bar Yc}  is not required anymore because the convex envelope of $\Phi$ is continuous, see Remark \ref{Remark convexity payoff1212} below.}
\end{Remark}

\subsection{Representation as a  BSDE with strong terminal condition}\label{subsec: attainability}

The previous section raises in particular one natural question: Does there exist an admissible control $\hat \alpha$ on the whole time interval $[0,T]$ allowing to match all time $t$-values of the minimal solution of a BSDE with weak terminal condition? Rephrasing, we wonder about the existence of a control $\hat\alpha$ in $\Ab_0$ such that 
\b*
\Yc^{\hat\alpha}_t &=& \Ec^g_{t}\left[\Phi(M^{\hat \alpha}_{T})\right] \;,  \qquad 0\le t\le T\;.
\e*
 Hereby, solving the BSDE with weak terminal condition $(\Psi,m_o,0)$ boils down to solving the classical BSDE with the optimal strong terminal one $\Phi(M^{\hat\alpha}_T)$: along the optimal path $\hat\alpha$, the compensator $\Kc^{\hat \alpha}$ of the   BSDE \reff{eq : BSDE de la cond minimale - BSDE}  must degenerate  to $0$. 

 {Not surprisingly, the existence of an optimal control requires the addition of convexity assumptions on the coefficients of the BSDE. We shall therefore assume that: \\

\textbf{ \Hgconv} \label{assumption:gconvex} For all $ (\lambda,m_{1},m_{2},t,y_1,y_2,z_1,z_2)\in [0,1]\times [0,1]^{2}\times[0,T]\times\R^2\times [\R^d]^2$, the following holds $\Pas$:
\b*
\Phi(\lambda m_{1}+(1-\lambda)m_{2})&\le &\lambda \Phi( m_{1} )+ (1-\lambda)\Phi( m_{2} )  \\
 g(t,\lambda y_1+(1-\lambda)y_2,\lambda z_1+(1-\lambda)z_2)&\leq &\lambda g(t,y_1,z_1) + (1-\lambda) g(t,y_2,z_2)  
 \e*

\begin{Remark} \label{RemConvexityEg}
{\rm We recall the following result which is based on standard comparison arguments, see e.g.  \cite[Proposition 7]{Rosazza-Gianin}: For any $\tau \in\Tc$, the map $\Ec^g_{\tau}[\Phi(\cdot)]:\Lb_0([0,1]) \to \Lb_0$ is convex under Assumption \textbf{\Hgconv}.}
\end{Remark}

\begin{Proposition}\label{prop : existence alpha optimal} Assume that Assumptions \textbf{\Hgconv} and \eqref{eq:hypcontiPhi} hold. Then, for any $(\tau, \alpha)\in \Tc  \x \Hb_{2}$, there exists $\hat \alpha^{\tau,\alpha}\in \Ab_{\tau}^{\alpha}$ such that 
$$
\Yc_{\tau}^{\alpha}=\Ec^g_{\tau}\left[\Phi(M^{\hat \alpha^{\tau,\alpha}}_{T})\right]=\Ec^g_{\tau,\tau'}\left[\Yc_{\tau'}^{{\hat \alpha}^{\tau,\alpha}}\right],\;\forall\; \tau'\in \Tc_{\tau}.
$$
\end{Proposition}

\begin{Remark} 
{\rm As detailed in Remark \ref{Remark convexity payoff} below, the convexity assumption on the terminal map $\Phi$ can be avoided in some cases. In particular, if  $\Phi$ is deterministic then it can  be replaced by its convex envelope. Then, only the convexity assumption on $g$ has to hold.}
\end{Remark}

\proof
Lemma \ref{lemma repres Y} below provides a sequence $(\alpha^n)_n$ valued in $\Ab_{\tau}^{\alpha}$ such that 
\begin{equation}
\label{eq:temp1_prop : existence alpha optimal}
\Yc_{\tau}^{\alpha}=\lim_{n\to \infty} \downarrow \Ec^g_{\tau}\left[\Phi(M^{\alpha^n}_{{T}})\right], \; \Pas 
\end{equation}
 
Since the sequence  $ (M^{\alpha^n}_{T})_{n}$ is bounded in $[0,1]$, we can   find  
sequences of non-negative real numbers $(\lambda_i^n)_{i\geq n}$ with $\sum_{i\geq n} \lambda_i^n=1$,  such that only a finite number of $\lambda_i^n$ do not vanish, for each $n$, and such that the sequence of convex combinations $(\tilde{M}^n_{T})_n$  given by 
 \be\label{temp123}
  \tilde{M}^n_{T}  &:=&\sum_{i\geq n} \lambda_i^n M^{\alpha^i}_{T}  
  \ee
  converges $\Pas$   to some $\hat M_{T}\in \Lb_{0}([0,1])$. By dominated convergence, the convergence holds in $\Lb_{2}$, in particular $E_{\tau}[\hat M_{T}]=M^{\alpha}_{\tau}$, and the martingale representation Theorem implies that we can find $\hat \alpha \in \Ab_{\tau}^{\alpha}$ such that $\hat M_{T}=M^{\hat \alpha}_{T}$.   Using the convexity of $\Phi$ and $g$, see Remark \ref{RemConvexityEg}, we deduce that 
$$
\tilde{Y}^n_{\tau}  : =\sum_{i\geq n} \lambda_i^n \Ec^{g}_{\tau}\left[\Phi(M^{  {\alpha}^i }_{T}) \right] \ge    \Ec^{g}_{\tau}\left[\Phi(\tilde M^{n}_{T}) \right].
$$
By \reff{eq:temp1_prop : existence alpha optimal}, $\tilde{Y}^n_{\tau}\to \Yc^{\alpha}_{\tau}$ $\Pas$ On the other hand, the convergence $\tilde M^{n}_{T}\to M^{\hat \alpha}_{T}$ in $\Lb_{2}$ combined with the boundedness and a.s. continuity of $\Phi$ implies that $\Phi(\tilde M^{n}_{T}) \to \Phi(M^{ \hat{\alpha} }_{T})$ in $\Lb_{2}$, after possibly passing to a subsequence. Therefore the convergence $\Ec^{g}_{\tau}\left[\Phi(\tilde M^{n}_{T}) \right] \to \Ec^{g}_{\tau}\left[\Phi(M^{ \hat{\alpha} }_{T}) \right]$ $\Pas$ follows by Proposition \ref{PropStability} below. This gives $\Yc^{\alpha}_{\tau}\ge \Ec^{g}_{\tau}\left[\Phi(M^{ \hat{\alpha} }_{T}) \right]$, while the converse holds by definition of $\Yc^{\alpha}_{\tau}$.

It remains to show that $\Yc_{\tau}^{\alpha}=\Ec^g_{\tau,\tau'}\left[\Yc_{\tau'}^{\hat \alpha}\right]$, for   $\tau'\in \Tc_{\tau}$. To see this, first note that the above implies that 
$\Yc^{\alpha}_{\tau}= \Ec^{g}_{\tau,\tau'}\left[\Ec^{g}_{\tau'}[\Phi(M^{ \hat{\alpha} }_{T}) ]\right]\ge \Ec^{g}_{\tau,\tau'}\left[\Yc_{\tau'}^{\hat \alpha}\right]$ by standard comparison arguments and the fact that $\Ec^{g}_{\tau'}[\Phi(M^{ \hat{\alpha} }_{T}) ]\ge \Yc_{\tau'}^{\hat \alpha}$ by definition. 
As above, we can find a sequence $(\hat{\alpha}^{n})\in \Ab_{\tau'}^{\hat \alpha}$ such that $\Ec^{g}_{\tau'}\left[\Phi(M^{\hat{\alpha}^{n}}_{T})\right]\to \Yc_{\tau'}^{\hat \alpha}$ $\Pas$ In view of Remark \ref{rem: borne unif Yc}, the convergence holds in $\Lb_{2}$ and Proposition \ref{PropStability} below implies 
$$
\Yc^{\alpha}_{\tau} \le \Ec^{g}_{\tau,\tau'}\left[\Ec^{g}_{\tau'}\left[\Phi(M^{\hat{\alpha}^n}_{T})\right]\right]\to \Ec^{g}_{\tau,\tau'}\left[\Yc_{\tau'}^{\hat \alpha}\right],
$$
where we used the fact that $\hat{\alpha}^{n}\in \Ab_{\tau'}^{\hat \alpha}\subset \Ab_{\tau}^{\alpha}$ to obtain the left hand-side. 
\ep 

\section{Main properties of the minimal initial condition process}
\label{SEC: Properties}
In this section, we emphasize remarkable properties of the map $\Yc_t\,:\,\mu\in\Lb_0([0,1],\Fc_t)\mapsto \Yc_t(\mu)$, for $t\in[0,T)$. We first derive the continuity of this map under a weak continuity assumption on $\Ec^g[\Phi(\cdot)]$. Then, we verify that this map (or more precisely   its l.s.c. envelope) is convex, and discuss the propagation of the convexity property to the time boundary $T-$. Finally, we retrieve, in this non-Markovian setting, a dual representation of the map $\Yc_0$, using solely probabilistic arguments.

\subsection{Continuity}\label{subsec: continuity}
\def\Rc{{\cal R}}

Our  continuity result is stated in terms of the quantities
$$
Err_t(\eta)  :=  \esssup \left\{  \Rc_{t}(M,M')\;:\; M,M'\in \Lb_0([0,1]) \;,\;\; E_{t}[|M-M'|^{2}]  \le \eta  \right\} , 
$$
defined for $\eta \in \Lb_0([0,1])$, in which  
$$
 \Rc_{t}(M,M'):=|\Ec_t^g[\Phi(M)] - \Ec_t^g[\Phi(M')] |.
$$
Observe that classical a priori estimates on BSDEs ensure that $Err_t(\eta_{n})\to 0$ as $\eta_{n}\to 0$ $\Pas$ with $(\eta_{n})_{n} \subset \Lb_0([0,1])$,  whenever $ \Phi $ is a deterministic Lipschitz map, see e.g. Proposition \ref{PropStability} below. This observation remains valid when  $ \Phi$ is simply continuous, via a classical convolution density argument for Lipschitz maps on bounded domains. The next result indicates that this property ensures the regularity of the map: $\mu \mapsto \Yc_t(\mu)$.

 \begin{Proposition}\label{prop: bar Yc est continue en P} 
Let $t<T$,  $\mu_{1},\mu_{2}\in \Lb_0([0,1],\Fc_{t})$. Then, 
$$
|\Yc_{t}(\mu_{1})-\Yc_{t}(\mu_{2})| \le Err_{t}(\Delta(\mu_{1},\mu_{2}))+ Err_{t}(\Delta(\mu_{2},\mu_{1})),
$$
where 
\b*
\Delta(\mu_{i},\mu_{j}):=(1-\frac{\mu_{i}}{\mu_{j}})  \1_{\{\mu_{i}<\mu_j\}}
+\frac{\mu_{i}-\mu_{j}}{1-\mu_j} \1_{\{\mu_{i}>\mu_j\}},\;i,j=1,2.
\e*
Moreover, 
\b*
|\Yc_{t}(\mu_{1})-\Yc_{t}(\mu_{2})|\1_{\{\mu_{1}=0\}} \le \Rc_{t}(\mu_{2},0)
\e*
and 
\begin{align*}
&|\Yc_{t}(\mu_{1})-\Yc_{t}(\mu_{2})|\1_{\{\mu_{1}=1\}} 
\\
&\le   \esssup \left\{ \Rc_t(1,M)\;:\; M\in \Lb_0([0,1]) \;,\;\; E_{t}[|1-M|^{2}]  \le 1-\mu_{2}  \right\}.
\end{align*}
In particular, if  $Err_t(\eta_{n})\to 0$ $\Pas$ as $\eta_{n}\to 0$ $\Pas$, for all $(\eta_{n})_{n} \subset \Lb_0([0,1])$, then   $\mu\in   \Lb_0((0,1),\Fc_{t}) \mapsto \Yc_{t}(\mu)$ is  continuous  for the sequential $\Pas$ convergence and the strong $\Lb_{2}$ convergence.
\end{Proposition}

\proof  {\bf Step 1.} Fix $\mu_{1},\mu_{2}\in \Lb_0([0,1],\Fc_{t})$.
Given  $\alpha_2\in\Ab_{t,\mu_2}$, we define
\b*
 \lambda &:=&  \frac{1-\mu_{1}}{1-\mu_2} \1_{\{\mu_2<\mu_{1}\}} +  \frac{\mu_{1}}{\mu_2}  \1_{\{\mu_{1} < \mu_2\}} +  \1_{\{\mu_{1} = \mu_2\}}\;, 
\e*
which is by construction valued in $[0,1]$. Since $M^{(t,\mu_2),\alpha_2}$ takes values in $[0,1]$,
\b*
M^{(t,\mu_{1}),\lambda\alpha_2}  =  \mu_{1}-\lambda\mu_2 + \lambda M^{(t,\mu_2),\alpha_2} \;\in\; [\mu_{1}-\lambda\mu_2,\mu_{1}+\lambda (1-\mu_{2})] \;\subset\;  [0,1] \;. 
\e*
In particular, $\lambda\alpha_2\in{\Ab}_{t,\mu_{1}}$. 
 Thus,  \reff{eq : def bar Y} leads to
\beq\label{temp112}
 \Yc_{t}(\mu_{1})
 \le
  \Ec^g_t[\Phi(M^{(t,\mu_2),\alpha_2}_T)] 
  + (   \Ec^g_t[\Phi(M^{(t,\mu_{1}),\lambda\alpha_2}_T)]  -   \Ec^g_t[\Phi(M^{(t,\mu_2),\alpha_2}_T)] )\;.
\eeq
Besides, 
$$
M^{(t,\mu_{1}),\lambda\alpha_2}_T  -  M^{(t,\mu_2),\alpha_2}_T
= \mu_{1} -\lambda\mu_2 +(\lambda -1)M^{(t,\mu_2),\alpha_2}_T
$$
so that, since $M^{(t,\mu_2),\alpha_2}_T$ belongs to $[0,1]$, we have
$$
\mu_{1}-1+\lambda(1-\mu_2) \leq M^{(t,\mu_{1}),\lambda\alpha_2}_T  -  M^{(t,\mu_2),\alpha_2}_T \le \mu_{1} -\lambda\mu_2.
$$
In addition,
\begin{align*}
\mu_{1} -\lambda\mu_2 \;=\; 0 \;, \qquad &\mbox{ if } \mu_{1}\;{<}\;\mu_2 \;, \mbox{ and }\\
\mu_{1}-1+\lambda(1-\mu_2) \;=\;0 \;, \qquad &\mbox{ if } \mu_{1}\;{\ge}\;\mu_2 \;.
\end{align*}
This directly leads to 
\b*
E_t[|M^{(t,\mu_{1}),\lambda\alpha_2}_T  -  M^{(t,\mu_2),\alpha_2}_T|] 
&\le& 
\Delta(\mu_{1},\mu_{2})
\;.
\e*
Since these two processes belong to $[0,1]$, we get 
$$
E_{t}[|M^{(t,\mu_{1}),\lambda\alpha_2}_T  -  M^{(t,\mu_2),\alpha_2}_T|^{2}] \le   \Delta(\mu_{1},\mu_{2}) .
$$
Hence, the arbitrariness of $\alpha_2\in{\Ab}_{t,\mu_2}$ together with  \reff{eq : def bar Y}  and  \reff{temp112} provides
\b*
 \Yc_{t}(\mu_{1})
 &\le&
 \Yc_{t}(\mu_2)
  + Err_t( \Delta(\mu_{1},\mu_{2}) ) \;.
\e*
 Interchanging the roles of $\mu_{1}$ and $\mu_2$ leads to 
 \b*
 \Yc_{t}(\mu_{2})
 &\le&
 \Yc_{t}(\mu_1)
  + Err_t( \Delta(\mu_{2},\mu_{1}) ) \;.
\e*

\noindent {\bf Step 2. } We next consider the case where $\Pro{\mu_{1}=0}>0$. Without loss of generality, we can assume that $\mu_{1}\equiv 0$. Fix $\alpha\in \Ab_{t,\mu_{2}}$. Since $\Ab_{t,\mu_{1}}=\{0\}$, $M^{(t,\mu_{2}),\alpha}_{T}\ge 0$ and $\Phi$ is non-decreasing, comparison implies that 
\b*
\Yc_{t}(0)=\Ec^{g}_{t}[\Phi(0)]\le \Ec^{g}_{t}[\Phi(M^{(t,\mu_{2}),\alpha}_{T})].
\e* 
In particular, $\Yc_{t}(0)=\Ec^{g}_{t}[\Phi(0)]\le \Yc_{t}(\mu_{2})\le \Ec^{g}_{t}[\Phi(M^{(t,\mu_{2}),0}_{T})]=\Ec_{t}^{g}(\Phi(\mu_{2}))$.\\

\noindent {\bf Step 3. } We  now consider the case where $\Pro{\mu_{1}=1}>0$. Again, we can assume that   $\mu_{1}\equiv 1$ so that $\Ab_{t,\mu_{1}}=\{0\}$. 
By comparison as above, one has 
$$
\Yc_{t}(1)=\Ec^{g}_{t}[\Phi(1)]\ge \Yc_{t}(\mu_{2}).
$$
On the other hand, since $M^{(t,\mu_{2}),\alpha}$ is a martingale taking values in $[0,1]$, we have
$$ E_t[|1-M_T^{(t,\mu_{2}),\alpha}|^2] \leq E_t[1-M_T^{(t,\mu_{2}),\alpha}]=1-\mu_2, \quad \alpha \in \Ab_{t,\mu_{2}},$$
from which the result follows.
\ep

\subsection{Convexity}\label{subsec: convexity}

In \cite{BoElTo09} and \cite{Mo10}, it is shown that the map $m\in [0,1]\mapsto \Yc_{0}(m)$ is convex. This is done in a Markovian framework using PDE arguments. In this section, we provide a probabilistic proof of this result which hereby extends to our setting. The result is stated for the lower-semicontinuous envelope $\Yc_{t*}$ of $\Yc_{t}$ defined as 
\be\label{eq: def Yct*}
\Yc_{t*}(\mu):=\lim_{\eps\to 0} \essinf\{\Yc_{t}(\mu'):~|\mu'-\mu|\le \eps,\;\mu'\in \Lb_{0}([0,1],\Fc_{t})\},
\ee
for any $t\in[0,T]$. We refer to Proposition \ref{prop: bar Yc est continue en P}, the discussion before it and to (ii) of Remark \ref{rem: after theorem charact bsde} for conditions ensuring that $\Yc_{*}=\Yc$. 
 
We first  make precise the notion of convexity adapted to our non-Markovian setting. Fix a time $t\in[0,T]$. 

\begin{definition}[$\Fc_t$-convexity] $\;$

\begin{enumerate}[(i)]
\item In the following, we say that a  {subset} $D\subset \Lb_{\infty}(\R,\Fc_{t})$ is  {$\Fc_{t}$-convex} if $\lambda \mu_{1}+(1-\lambda)\mu_{2}\in D$, for all $\mu_{1},\mu_{2}\in D$ and $\lambda\in \Lb_0([0,1],\Fc_{t})$.

\item Let $D$ be an $\Fc_{t}$-convex subset of $\Lb_{\infty}(\R,\Fc_{t})$. 
A  {map} $\Jc : D\mapsto \Lb_2(\R,\Fc_{t})$ is said to be  {$\Fc_{t}$-convex} if 
$$
{\rm Epi}(\Jc):=\left\{(\mu,Y)\in D\x \Lb_2(\R,\Fc_{t}):~ Y\ge\Jc(\mu)\right\} 
$$
is $\Fc_{t}$-convex.

\item 
Let 
$ {\rm Epi}^{c}(\Yc_{t})$ be the set of elements of the form $\sum_{n\le N}  \lambda_{n}  (\mu_{n},Y_{n})$ with  $
 (\mu_{n},Y_{n},\lambda_{n})_{n\le N}\subset   {\rm Epi}(\Yc_{t})\x \Lb_0([0,1],\Fc_{t})$ such that $\sum_{n\le N}\lambda_{n}=1$, for some $N\geq 1$. We then denote by $\overline {\rm Epi}^{c}(\Yc_{t})$
 its closure in $\Lb_2$. Finally, the   {$\Fc_{t}$-convex envelope} of $\Yc_{t}$ is defined as 
 \be\label{eq: def Yctconv}
 \Yc_{t}^{c}(\mu):=\essinf\{Y\in \Lb_2(\R,\Fc_{t}): (\mu,Y)\in \overline {\rm Epi}^{c}(\Yc_{t})\}.
 \ee
 \end{enumerate}
\end{definition}
  
 We can now state the convexity property. It requires a right continuity property in time,  which holds under the conditions of Theorem \ref{thm : BSDE characterization bar Yc}{(ii)}, also recall  (ii) of Remark \ref{rem: after theorem charact bsde}. 
  
 \begin{Proposition}\label{prop: bar Yc est convex en P} 
Assume that $\Yc_{t}(\mu)=\Yc_{t+}(\mu)$  for any $\mu \in \Lb_{0}([0,1],\Fc_{t})$ and $t< T$. Then, the map $\mu\in \Lb_0([0,1],\Fc_{t}) \mapsto \Yc_{t*}(\mu)$ is $\Fc_{t}$-convex, for all $t<T$. 
\end{Proposition}

\proof Fix $t\in[0,T)$ and set $D:=\Lb_0([0,1],\Fc_{t})$ for ease of notations. The proof is divided in several steps.\\

\no {\bf Step 1.} {\sl  
$
(\mu,\Yc^{c}_{t}(\mu))\in \overline {\rm Epi}^{c}(\Yc_{t})
$, for all $\mu \in D$.   
}
\vs2

Indeed, the family $F:=\{Y\in \Lb_2(\R,\Fc_{t})$ $:$ $(\mu,Y)\in \overline {\rm Epi}^{c}(\Yc_{t})\}$ is directed downward (for every fixed element $\mu$ in $D$) since $Y^{1}\Indi{\{Y^{1}\le Y^{2}\}}+ Y^{2}\Indi{\{Y^{1}> Y^{2}\}} \in F$, by $\Fc_{t}$-convexity of $\overline {\rm Epi}^{c}(\Yc_{t})$, for all $Y^{1},Y^{2}\in F$.  It then follows from \cite[Proposition VI.1.1]{Ne75} that there exists a sequence $(Y^{n})_{n\ge 1} \subset F$ such that $Y^{n}\downarrow  \Yc^{c}_{t}(\mu)$ $\Pas$ Since $Y^{1}$ and $\Yc^{c}_{t}(\mu) \in \Lb_{2}$,  the monotone convergence Theorem implies that  $Y^{n}\to  \Yc^{c}_{t}(\mu)$  in $\Lb_2$, as $n$ goes to infinity. The set $\overline {\rm Epi}^{c}(\Yc_{t})$ being closed in $\Lb_2$, this proves our   claim. 
\vs2

\no {\bf Step 2.} {\sl  Let $\eta\in \Sb_{2}$ be as in Remark \ref{rem: borne unif Yc}. Then,   $|\Yc^{c}_{t}(\mu)|\le \eta_{t}$, for all $t\le T$ and $\mu \in D$. 
}
\vs2

We first observe that $\Yc\ge \Yc^{c}$ by construction.  Remark \ref{rem: borne unif Yc} thus implies that $ \Yc^{c}_{t}(\mu)\le \eta_{t}$. On the other hand, let $(Y^{n})_{n\ge 1}$ be as in the step above. We claim that it satisfies $Y^{n}\ge -\eta_{t}$, for each $n\ge 1$.  Then, the lower bound $ \Yc^{c}_{t}(\mu)\ge -\eta_{t}$ is obtained by passing to the limit. To see this, it suffices to prove this property for any $Y  \in \Lb_2(\R,\Fc_{t})$ such that  $(\mu,Y)\in \overline {\rm Epi}^{c}(\Yc_{t})$. But, such an element $(\mu,Y)$ is obtained by taking the $\Lb_{2}$ limit of elements of the form $\sum_{n\le N}  \lambda_{n}  (\mu_{n},Y_{n})$ with  $
 (\mu_{n},Y_{n},\lambda_{n})_{n\le N}\subset   {\rm Epi}(\Yc_{t})\x \Lb_0([0,1],\Fc_{t})$, such that $\sum_{n\le N}\lambda_{n}=1$. Each $Y_{n}$ of the latter family is bounded from below by $-\eta_{t}$ by Remark \ref{rem: borne unif Yc}, and hence so is $Y$.
\vs2
  
\no {\bf Step 3.} {\sl The map   $\mu\in D\mapsto \Yc^{c}_{t}(\mu)$ is $\Fc_{t}$-convex.}
\vs2

Fix $\mu^{1}, \mu^{2} \in D$ and $\lambda\in \Lb_0([0,1],\Fc_{t})$. Step 1 implies that $(\mu^{i},\Yc^{c}_{t}(\mu^{i}))\in \overline {\rm Epi}^{c}(\Yc_{t}) $ for $i=1,2$. Clearly,  $ \overline {\rm Epi}^{c}(\Yc_{t})$   is $\Fc_{t}$-convex. It follows that  $(\lambda \mu^{1}+(1-\lambda ) \mu^{2},$ $\lambda \Yc^{c}_{t}(\mu^{1})$ $+  (1-\lambda ) \Yc^{c}_{t}(\mu^{2}) )$ $\in$  $  \overline {\rm Epi}^{c}(\Yc_{t})$, so that 
$\lambda \Yc^{c}_{t}(\mu^{1})$ $+  (1-\lambda ) \Yc^{c}_{t}(\mu^{2}) \ge \Yc^{c}_{t}(\lambda \mu^{1}+(1-\lambda ) \mu^{2})$. Now, for any $Y^i$ such that $(\mu^{i},Y^{i})\in {\rm Epi}(\Yc^{c}_{t})$,  one has  $Y^{i}\ge    \Yc^{c}_{t}(\mu^{i})$, $i=1,2$. This fact combined with the previous inequality thus implies $\lambda Y^{1}$ $+  (1-\lambda ) Y^{2} \ge \Yc^{c}_{t}(\lambda \mu^{1}+(1-\lambda ) \mu^{2})$. This means that ${\rm Epi}(\Yc^{c}_{t})$ is $\Fc_{t}$-convex. 
\vs2
 
\no{\bf Step 4. } {\sl $\Yc_{t*}(\mu)\ge \Yc^{c}_{t}(\mu)$, for all $\mu\in D$.} 
\vs2
\def\Ic{{\cal I}}

Fix $\varepsilon>0$ and set $D_{\mu}^{\eps}:= \{ {\mu'\in \Lb_{0}([0,1],\Fc_{t}), \;|\mu'-\mu|\le \eps}\}$. It follows from Remark \ref{rem: conca sur event set} that the family 
$ \{\Yc_{t}(\mu'):~ \mu'\in D_{\mu}^{\eps}\} $ 
is directed downward. Then, we can find a sequence $(\mu_{n}^{\eps})_{n\ge 1} \subset D_{\mu}^{\eps}$ such that 
\b*
 \Yc_{t}(\mu^{\eps}_{n})\to Z_{\eps}(\mu):=\essinf\{\Yc_{t}(\mu'):~ \mu'\in D_{\mu}^{\eps}\}\quad\Pas
\e*
Since $(Z_{\eps}(\mu))_{\eps>0}$ is non-decreasing,   $\lim_{N\to \infty} Z_{1/N}(\mu)=\Yc_{t*}(\mu)$, recall \reff{eq: def Yct*}. Note that Remark \ref{rem: borne unif Yc} implies that $( \Yc_{t}(\mu^{1/N}_{n}))_{n\ge 1} \to_{n} Z_{1/N}(\mu)$ in $\Lb^{2}$ and define
$$
k_{N}:=\min\{n\ge 1: \| \Yc_{t}(\mu^{1/N}_{n})-Z_{1/N}(\mu)\|_{\Lb^{2}}\le 1/N\}.
$$
Then, $(\mu_{k_{N}}^{1/N}, \Yc_{t}(\mu^{1/N}_{k_{N}}))\to (\mu,\Yc_{t*}(\mu))$ in $\Lb^{2}$ as $N\to \infty$. Since $  {\rm Epi} (\Yc_{t})\subset \overline {\rm Epi}^{c}(\Yc_{t})$ and  the latter is closed under $\Lb^{2}$-convergence, this implies that $ (\mu,\Yc_{t*}(\mu)) \in \overline {\rm Epi}^{c}(\Yc_{t})$. We conclude by appealing to the definition of $\Yc_{t}^{c}$ in \reff{eq: def Yctconv}.

\vs2
\no{\bf Step 5. } {\sl $\Yc^{c}_{t}(\mu)\ge  \Yc_{t*}(\mu)$, for all $\mu\in D$.} 
\vs2

In view of Steps 3 and 4, the result of Step 5 actually proves that $\Yc_{t*}=\Yc^{c}_{t}$ is $\Fc_{t}$-convex. 

\noindent We now proceed to the proof of Step 5 which is itself divided in two parts. 

\no{\bf Step 5.a } It follows from Step 1, that there exists a sequence  
\be\label{eq: proof convexity def sequence Pn Yn Lambda}
(\mu_{n},Y_{n},\lambda^{N}_{n})_{n\ge 1,N\ge 1}\subset   {\rm Epi}(\Yc_{t})\x \Lb_0([0,1],\Fc_{t})
\ee
 such that $\sum_{n\le N}  \lambda^{N}_{n}=1$, for all $N$, and 
 \be\label{eq: proof convexity def sequence Pn Yn Lambda conv}
 (\hat \mu_{N},\hat Y_{N}):=\sum_{n\le N}  \lambda^{N}_{n}  (\mu_{n},Y_{n})\to (\mu, \Yc_{t} ^{c}(\mu)) \;\mbox{ in }\Lb_2.
 \ee
Fix $N\ge 1$ and $\eps>0$. Let $\hat \alpha^{N}\in \Hb_{2}$ be such that $\hat \mu_{N}=m_{o}+\int_{0}^{t} \hat \alpha^{N}_{s} dW_{s}$.     
 Since the family $(\lambda_{n}^{N})_{n\le N}$ is composed of $\Fc_{t}$-measurable random variables summing to $1$, one can find $\alpha^{N}\in \Hb_{2}$ and a random variable $\xi^{\eps}_{N}\in \Lb_2(\Fc_{t+\eps})$ such that 
\beq
\label{eq:combconvparprobacond}
 \hat \mu_{N} +\int_{t}^{t+\eps}\alpha^{N}_{s} dW_{s}= \xi^{\eps}_{N}
\;\mbox{ and }\;
  \Pro{\xi^{\eps}_{N}=\mu_{n}|\Fc_{t}}=\lambda^{N}_{n}\;, \;\;\mbox{ for }\;n\le N.
\eeq
 Without loss of generality, we can assume that $\alpha^{N}=\hat \alpha^{N}$ $dt\times d\P$ on $[0,t]$. 
Then, (i) of Theorem \ref{thm : BSDE characterization bar Yc} and Remark \ref{rem: conca sur event set} yield
\begin{align}
\Yc_{t}(\hat \mu_{N}) &=\Yc_{t}^{\hat \alpha^{N}}\le \Ec_{t,t+\eps}^{g}(\Yc^{\alpha^{N}}_{t+\eps})=\Ec_{t,t+\eps}^{g}(\Yc_{t+\eps}(\xi^{\eps}_{N})) \nonumber\\
&= \Ec_{t,t+\eps}^{g}\left(\sum_{n\le N} 1_{\xi^{\eps}_{N}=\mu_{n}}\Yc_{t+\eps}(\mu_{n})\right)\;.\label{eq: proof convexity serie inega avant eps to 0}
\end{align}
We claim that
\be\label{eq: proof convexity claim conv as eps to 0}
\liminf_{\eps\to 0}\Ec_{t,t+\eps}^{g}\left(\sum_{n\le N} 1_{\xi^{\eps}_{N}=\mu_{n}}\Yc_{t+\eps}(\mu_{n})\right)\le  \sum_{n\le N}  \lambda^{N}_{n}  \Yc_{t}(\mu_{n}).
\ee
Then,  \reff{eq: proof convexity serie inega avant eps to 0}, \reff{eq: proof convexity claim conv as eps to 0}, \reff{eq: proof convexity def sequence Pn Yn Lambda} and \reff{eq: proof convexity def sequence Pn Yn Lambda conv} lead to 
$$
\Yc_{t}(\hat \mu_{N}) \le \sum_{n\le N}  \lambda^{N}_{n}  \Yc_{t}(\mu_{n}) \le  \sum_{n\le N}  \lambda^{N}_{n} Y_{n}=\hat Y_{N}.
$$ 
Appealing to \reff{eq: proof convexity def sequence Pn Yn Lambda conv}, we deduce that
$$
  \liminf_{N\to \infty } \Yc_{t}(\hat \mu_{N}) \le \Yc_{t}^{c}(\mu).
$$ 
Since $\hat \mu_{N}\to \mu$ $\Pas$, this together with Remark \ref{rem: conca sur event set} implies that  
\begin{align*}
Z_{\eps}(\mu)&\le \liminf_{N\to \infty} \Yc_{t}(\bar \mu_{N})= \liminf_{N\to \infty} \left(\Yc_{t}(\hat \mu_{N})\1_{\{|\hat \mu_{N}-\mu|\le \eps\}} + \Yc_{t}(\mu)\1_{\{|\hat \mu_{N}-\mu|> \eps\}} \right)\\ &\le    \Yc_{t}^{c}(\mu),  
\end{align*}
for all $\eps>0$, where     
$$
\bar   \mu_{N}:= \hat \mu_{N} \1_{\{|\hat \mu_{N}-\mu|\le \eps\}} + \mu\1_{\{|\hat \mu_{N}-\mu|> \eps\}} \in D^{\eps}_{\mu}, 
$$
see Step 4 for the definitions of $Z_{\eps}(\mu)$ and $D^{\eps}_{\mu}$. 
 Since $Z_{\eps}(\mu)\uparrow \Yc_{t*}(\mu)$ as $\eps$ goes to $0$ by \eqref{eq: def Yct*}, this shows the required result. 
\vs2

\no{\bf Step 5.b } It finally remains to prove the claim \reff{eq: proof convexity claim conv as eps to 0}.
 
\noindent Remark \ref{rem: borne unif Yc} and (ii) of Proposition \ref{prop: majo bsde appendix} in the Appendix imply that  
\b*
 \Ec_{t,t+\eps}^{g}\left(\sum_{n\le N} 1_{\xi^{\eps}_{N}=\mu_{n}}\Yc_{t+\eps}(\mu_{n})\right)
 &\le& E_{t}\left[\sum_{n\le N} 1_{\xi^{\eps}_{N}=\mu_{n}}\Yc_{t+\eps}(\mu_{n}) \right]+  \eta_{\eps} \\
& \le& E_{t}\left[\sum_{n\le N} 1_{\xi^{\eps}_{N}=\mu_{n}}\Yc_{t}(\mu_{n}) \right]+  \eta_{\eps} \\
&+& \sum_{n\le N}  E_{t}\left[|\Yc_{t+\eps}(\mu_{n})-\Yc_{t}(\mu_{n})| \right],
 \e*
where $\eta_{\eps} \to 0$ $\Pas$ as $\eps\to 0$. The right-hand side of \eqref{eq:combconvparprobacond} then leads to 
\b*
 \Ec_{t,t+\eps}^{g}\left(\sum_{n\le N} 1_{\xi^{\eps}_{N}=\mu_{n}}\Yc_{t+\eps}(\mu_{n})\right)
& \le& \sum_{n\le N}  \lambda^{N}_{n}\Yc_{t}(\mu_{n})+  \eta_{\eps}  \\
&+& \sum_{n\le N}  E_{t}\left[|\Yc_{t+\eps}(\mu_{n})-\Yc_{t}(\mu_{n})| \right].
 \e*
 Recall  that    $\Yc_{t+}(\mu_{n})=\Yc_{t}(\mu_{n})$ by assumption,  and that $(\Yc(\mu_{n}))_n$  is bounded by some $\eta\in \Sb_{2}$, see Remark \ref{rem: borne unif Yc}. Sending $\eps\to 0$  in the above inequality  and appealing to the Lebesgue dominated convergence Theorem  proves \reff{eq: proof convexity claim conv as eps to 0}.
 \ep\\

 In the context of PDEs, convexity in the domain propagates up to the boundary, which leads to a boundary layer phenomenon. In  \cite{BoElTo09} and \cite{Mo10} this translates in the fact that the natural $T$-time boundary condition should be stated in terms of the $m$-convex envelope of $\Phi$. We observe hereafter that this property extends to our non-Markovian setting, whenever $\Phi$ is deterministic. 
 
We recall from Theorem \ref{thm : BSDE characterization bar Yc} (i) that $\Yc$ {is l\`adl\`ag on countable sets. Under the following condition, it will actually be c\`adl\`ag up to undistinguishability}. As opposed to Proposition \ref{prop: bar Yc est convex en P}, we shall not need to impose any right-continuity for the following.  

\begin{Proposition}\label{prop: regularisation au bord} 
Assume that $\Phi$ is deterministic and let  $\hat \Phi$ denote its convex envelope. Then, 
$$
{\lim_{t\uparrow T}\Yc^{\alpha}_{t}  } =\hat \Phi(M^{\alpha}_{T})
\;\mbox{ and }\;
\Yc_{\tau}^{\alpha}=\essinf_{\alpha'\in \Ab_{\tau}^{\alpha} }  \Ec^g_{\tau}\left[ \hat \Phi(M^{\alpha^{'}}_{T})\right],
$$
for all $\alpha \in \Ab_{0}$ and $\tau \in \Tc$ such that $\tau<T$.
\end{Proposition}

Before proving this result, let us make some observations. 

\begin{Remark}\label{Remark convexity payoff1212}
{\rm Since $\Phi$ is non-decreasing, its convex envelope is continuous on $[0,1)$. Moreover, $\Phi$ is left-continuous, so that $\hat \Phi$ has to be continuous at $1$ as well.}
\end{Remark}

\begin{Remark}\label{Remark convexity payoff}
{\rm In Section \ref{subsec: attainability}, we observed that the essential infimum in the dynamic programming principle is attained whenever $\Phi$ and $g$ are convex. Hence, the previous proposition allows straightforwardly to avoid the convexity requirement on $\Phi$, whenever it is deterministic.}
\end{Remark}

\begin{Remark}\label{Remark quantile hedging}
{\rm The proof below can easily be adapted to the case where $\Phi(\omega,m)=\phi(m)\xi(\omega)$ for some non-negative random variable $\xi$ and a deterministic map $\phi$.  This is due to the fact that the $m$-convex envelope of $\Phi$ is fully characterized by the convex envelope $\hat \phi$ of $\phi$: $\hat \Phi(\omega,m)=\hat \phi(m)\xi(\omega)$. This allows one to follow the construction used in our proof.  In particular, in the quantile hedging problem of F\"{o}lmer and Leukert \cite{FoLe99}, one has $\Phi(\omega,m)=\1_{\{m>0\}}\xi(\omega)$ ($m\in [0,1]$), with $\xi$ taking non-negative values, so that $\hat \Phi(\omega,m)=m\xi(\omega)$,  see also \cite{BoElTo09}. 
}
\end{Remark}

\noindent {\bf Proof of Proposition \ref{prop: regularisation au bord}.} We prove each assertion separately. 

\no{\bf Step 1.} By definition of the convex envelope, we can find a measurable map $m\in [0,1]\mapsto (\underline \wp(m),\overline \wp(m),\eps(m))\in [0,1]^{3}$ such that $\underline \wp(m)\le m \le \overline\wp(m)$,  $\eps(m) \underline \wp(m)+(1-\eps(m))\overline \wp(m)=m$ and 
$$
\hat \Phi(m)=\eps(m)\Phi(\underline \wp(m))+(1-\eps(m))\Phi(\overline \wp(m))\;,
$$ 
for any $m\in[0,1]$.
Let $t_{n}\uparrow T$. Then, one can find $\alpha^{n}\in \Ab_{t_{n}}^{\alpha}$ and $\xi^{n}\in \Lb_0([0,1])$ such that $M^{\alpha^{n}}_{T}=M^{\alpha}_{t_{n}}+\int_{t_{n}}^{T} \alpha^{n}_{s}dW_{s}=\xi^{n}$, where $\Pro{\xi^{n}=\underline \wp(M^{\alpha}_{t_{n}})|\Fc_{t_{n}}}$ $=$ $\eps(M^{\alpha}_{t_{n}})$ and 
$\Pro{\xi^{n}=\overline \wp(M^{\alpha}_{t_{n}})|\Fc_{t_{n}}}=1-\eps(M^{\alpha}_{t_{n}})$. 
It follows from the above and (iii) of Proposition \ref{prop: majo bsde appendix} in the Appendix that 
\b*
\Yc^{\alpha}_{t_{n}}&\le& E_{t_{n}}\left[\Phi(\xi^{n}) \right]+\eta_{n}
= \hat\Phi( M^{\alpha}_{t_{n}})+\eta_{n},
\e*
where $\eta_{n}\to 0$ as $n\to \infty$.
Since $\Yc$ {is l\`adl\`ag on countable sets} (by Proposition \ref{prop : aggregation cadlag}), passing to the limit  implies that 
 \be\label{temp1245}
{\lim_{n\to \infty}\Yc^{\alpha}_{t_{n}}}&\le& \hat \Phi(M^{\alpha}_{T}).
 \ee
We now prove the converse inequality. 
We use (iii) in Proposition \ref{prop: majo bsde appendix} in the Appendix  and   Jensen's inequality to deduce that 
\b*
Y^{\alpha'}_{t_{n}}&:=&\Eg{t_{n},T}{ \Phi(M^{\alpha'}_{T})}\ge  E_{t_{n}}\left[\hat \Phi(M^{\alpha'}_{T}) \right]-\bar\eta_{n}
\ge    \hat \Phi(M^{\alpha}_{t_{n}})-\bar\eta_{n}\;, \quad \alpha'\in \Ab_{t_{n}}^{\alpha}, 
\e*
where $\bar\eta_{n}\to 0$ as $n\to \infty$. Combining the arbitrariness of $\alpha'\in \Ab_{t_{n}}^{\alpha}$ with the  l\`adl\`ag property of $\Yc$ {on countable sets}, we get that 
$$
 {\lim_{n\to \infty}\Yc^{\alpha}_{t_{n}}}\ge \liminf_{n\to \infty } \essinf_{\alpha'\in \Ab_{t_{n}}^{\alpha}} Y^{\alpha'}_{t_{n}} \ge \hat \Phi(M^{\alpha}_{T})\;.
$$
 
\no{\bf Step 2.} It follows from Theorem \ref{thm : BSDE characterization bar Yc} (i)  that 
$$
\Yc_{\tau}^{\alpha}=\essinf_{\alpha'\in \Ab_{\tau}^{\alpha} } \Eg{\tau,t_{n}\vee \tau }{ \Yc_{t_{n}\vee \tau}^{\alpha'}}\;, \quad n\in\N\;.
$$
The process $\Yc_{.\vee \tau}^{\alpha'}$ being l\`adl\`ag {on the set $\{t_{n},n\ge 1\}$},  {$\lim_{n\to \infty}\Yc_{t_{n}\vee \tau}^{\alpha'}$}  is well-defined {and coincides with $\lim_{n\to \infty}\Yc^{\alpha'}_{t_{n}}$}. 
Moreover, it follows from the bound in Remark \ref{rem: borne unif Yc} that the convergence holds in $\Lb_2$. In view of the stability result of Proposition \ref{PropStability} and Step 1. above,  passing to the limit as $n\to \infty$ leads to 
$$
\Yc_{\tau}^{\alpha}\le \essinf_{\alpha'\in \Ab_{\tau}^{\alpha} } \Eg{\tau} {{ \lim_{n\to \infty}\Yc^{\alpha'}_{t_{n}}}}= \essinf_{\alpha'\in \Ab_{\tau}^{\alpha} } \Eg{\tau}{\hat \Phi(M^{\alpha'}_{T})}.
$$
Since $\Phi\ge \hat \Phi$, the reverse inequality holds by definition of $\Yc_{\tau}^{\alpha}$ {in \reff{eq : def bar Y}}. {Since $\hat \Phi$ is continuous by Remark \ref{Remark convexity payoff1212}, we can now appeal to the second statement of Proposition   \ref{prop : aggregation cadlag} to assert that, up to indistinguishability, $\Yc^{\alpha}$ is c\`adl\`ag, so that 
$\lim_{t\uparrow T} \Yc^{\alpha}_{t}=\lim_{n\to \infty}\Yc^{\alpha}_{t_{n}}$.}
\ep

\subsection{Dual representation}

In this section, we provide a dual formulation for the minimal initial condition  at time $0$,  $m\mapsto \Yc_{0}(m)$. It requires the introduction of the Fenchel transforms of  $g$ and $  \Phi$.

We therefore define  
$$
\tilde \Phi~:~ (\omega,l)\in \Omega\x \R \mapsto \sup_{m\in [0,1]} (ml-\Phi(\omega,m))
$$
and 
$$
\tilde g  ~:~ (\omega,t,u,v)\in \Omega\x [0,T]\x \R\x \R^{d} \mapsto\sup_{(y,z)\in \R\x \R^{d}} \left(yu + z^{\top}v-g(\omega,t,y,z)\right). 
$$
\begin{Remark}\label{rem: domaine tilde Phi et tilde g}
{\rm It follows from the assumption \Hg~that the domain of $\tilde g(\omega,t,\cdot)$, $\dom{\tilde g(\omega,t,\cdot)}$, is contained in $[-K_{g},K_{g}]^{d+1}$   for $\Pae$ $\omega \in \Omega$ and all $t\le T$. The assumption \HPsi~ensures that the domain of $\tilde \Phi(\omega,\cdot)$ is the all real line, $\Pas$.}
\end{Remark}
 
 In the following, we denote by $\Lambdab$ the set of predictable processes $\lambda$ with values  in $\R\x \R^{d}$ such that 
  $\lambda_{t}(\omega) \in \dom{\tilde g(\omega,t,\cdot)}$ for $\Leb \x  \P$-a.e. $(\omega,t)\in \Omega\x [0,T]$. 
  
  To $\lambda=(\nu,\vc)\in \Lambdab$, we associate the process $L^{\lambda}$ defined by 
$$
L_{t}^{\lambda} =1+\int_{0}^{t}L_{s}^{\lambda}\nu_{s}ds+\int_{0}^{t}L_{s}^{\lambda}\vc_{s}dW_{s}\;,\;t\in [0,T].
$$ 

Our  dual formulation for $\Yc_{0}$ is stated in terms of 
$$
    \Xc_{0}(l):=\inf_{\lambda\in \Lambdab} {\rm X}_{0}^{l,\lambda},\;l>0,
$$
where 
$$
{\rm X}^{l,\lambda}_{0} :=E\left[ \int_{0}^{T} L_{s}^{\lambda} \tilde g(s,\lambda_{s}) ds + L_{T}^{\lambda} \tilde \Phi(l/L_{T}^{\lambda})\right] \;, \quad \lambda\in\Lambdab\;,\quad l>0.
$$

The fact that the Fenchel transform of $\Xc_{0}$ provides a lower bound for $\Yc_{0}$ is straightforward, and detailed in   Proposition \ref{prop: borne inf dualite}  below for the convenience of the reader.   For ease of notations, we now write $\Ab_{m}$ for $\Ab_{0,m}$, $M^{m,\alpha}$ for $M^{(0,m),\alpha}$, and denote by  $(Y^{m,\alpha},Z^{m,\alpha})$   the solution of the $\BSDE(g,\Phi(M^{m,\alpha}_{T}))$, $\alpha\in \Ab_{m}$.

\begin{Proposition}\label{prop: borne inf dualite} 
$\Yc_{0}(m)\ge \sup_{l>0} \left( lm -     \Xc_{0}(l)\right)$, for all $m\in [0,1]$.
\end{Proposition}

\proof  Fix $\alpha\in \Ab_{m}$ and $\lambda=(\nu,\vc)\in \Lambdab$.    Then, it follows from the definition of $\tilde \Phi$ and $\tilde g$ that
\begin{align*}
E\left[Y^{m,\alpha}_{T}L_{T}^{\lambda}\right]&=
Y^{m,\alpha}_{0} \\& +E\left[ \int_{0}^{T} L_{s}^{\lambda} \left(  \nu_{s} Y^{m,\alpha}_{s} +  \vc_{s}^{\top} Z^{m,\alpha}_{s} - g(s,Y^{m,\alpha}_{s},Z^{m,\alpha}_{s})\right)ds\right]
\\
&\le
Y^{m,\alpha}_{0} +E\left[ \int_{0}^{T} L_{s}^{\lambda} \tilde g(s,\lambda_{s}) ds \right],
\end{align*}
and   
\b*
Y^{m,\alpha}_{T}L_{T}^{\lambda}=\Phi(M^{m,\alpha}_{T})L_{T}^{\lambda}\ge  l   M^{ m,\alpha}_{T}-L_{T}^{\lambda}\tilde \Phi(l/L_{T}^{\lambda}),
\e*
for $l>0$. Note that, in the above, we have cancelled the expectation of the local martingale part  $\int_0^T(L^{\lambda}_s Z^{m,\alpha}_s+Y^{m,\alpha}_sL^{\lambda}_s\vartheta_s)dW_s$ although $L^{\lambda}  Z^{m,\alpha}$ might not belong to $\Hb_2$. If not, one may use a localization argument since all other terms belongs to $L^1$ uniformly in time.
Combining the above {and using the martingale property of $M^{ m,\alpha}$} yields
$$
Y^{m,\alpha}_{0}\ge lm - E\left[ \int_{0}^{T}  L_{s}^{\lambda} \tilde g(s,\lambda_{s}) ds +L_{T}^{\lambda} \tilde \Phi(l/L_{T}^{\lambda})\right]
=lm-X_{0}^{l,\lambda} \;.
$$
The result follows from the arbitrariness of $l>0$, $\lambda\in \Lambdab$, and $\alpha \in \Ab_{m}$.
\ep
\vs2

We now show that equality is satisfied in Proposition  \ref{prop: borne inf dualite} whenever existence holds in the dual problem. This is proved under the following assumptions. Let $C^{1}_{b}$ be the set of continuously differentiable maps with bounded first derivatives.\\

\def\Hd#1{${\rm \bf (H_{d}^{#1})}$}

{\bf Assumption \Hd1}  The following holds for 
$\Leb \x \P$-a.e.  $(t,\omega) \in [0,T]\x \Omega$:
\begin{enumerate}[(a)]
\item\label{item: Hd1 C1b} the maps $ \tilde \Phi(\omega,\cdot )$ and $\tilde g(\omega,\cdot)$ are $C^{1}_{b}$ on their domain,  and $\dom{\tilde g(\omega,t,\cdot)} $ is closed;
\item\label{item: Hd1 borne nabla}  $|\nabla  \tilde \Phi(\omega,\cdot)| + |\nabla \tilde g(\omega,t,\cdot)|\le \chi_{\tilde \Phi,\tilde g}(\omega)$, for some $\chi_{\tilde \Phi,\tilde g}\in \Lb_{2}(\R)$;
\item\label{item: Hd1 bi fenchel egale Phi} $\Phi(\omega,m)=\sup\limits_{l>0} \left(lm-\tilde \Phi(\omega,l)\right)$, for all $m\in [0,1]$;
\item\label{item: Hd1 bi fenchel egale g}  $ g(\omega,t,y,z)=\max\limits_{(u,v)\in \dom{\tilde g(\omega,t,\cdot)}} \left(yu+z^{\top} v - \tilde g(\omega,t,u,v) \right)$, for all $(y,z)\in \R\x \R^{d}$.
\end{enumerate}
\vs2
In the above, $\nabla  \tilde \Phi$ and $\nabla \tilde g$ stands for the gradient with respect to $l$ and $(u,v)$ respectively.  

Note that \reff{item: Hd1 C1b} and \reff{item: Hd1 borne nabla} are of technical nature, while \reff{item: Hd1 bi fenchel egale Phi} and \reff{item: Hd1 bi fenchel egale g} mean that $\Phi$ and $g$ are  convex, i.e. coincide with their bi-dual. The latter is a minimal requirement if one wants the duality to hold. 

\begin{Proposition}\label{prop dualite existence duale} Let Assumption \Hd1 hold. Assume further that there exists $\hat l>0$ and $\hat \lambda \in \Lambdab$  such that 
\be
&  \sup\limits_{l>0} \left(lm -    \Xc_{0}(l)\right)= \hat lm-   \Xc_{0}(\hat l)= \hat lm- {\rm X}_{0}^{\hat l,\hat\lambda}.&
 \ee
 Then,  there exists $\hat \alpha\in \Ab_{m}$ such that 
 $$
  \Yc_{0}(m)=Y^{m,\hat \alpha}_{0}=  \hat lm-    \Xc_{0}(\hat l). 
 $$ 
It satisfies
 \begin{equation}\label{eq: prop dualite existence duale - relation variables}
  g(\cdot,Y^{m,\hat \alpha},Z^{m,\hat \alpha})=  \hat \lambda ^{\top}(Y^{m,\hat \alpha},Z^{m,\hat \alpha})-\tilde g (\cdot,\hat \lambda)  
\mbox{ , } \Phi(M^{m,\hat \alpha}_{T})=\frac{M^{m,\hat \alpha}_{T}\hat l}L^{\hat \lambda}_{T}-  \tilde \Phi(\hat l/L^{\hat \lambda}_{T}).
 \end{equation}
 \end{Proposition}

Before to provide the proof, let us make the following observation which pertains for the case of a linear driver $g$.

\begin{Remark}{\rm  Assume that $g$ is linear, i.e. there exist bounded predictable processes $A^{Y}$ and $A^{Z}$ such that
$g: (\omega,t,y,z) \mapsto g(\omega,t,0,0)+ A^{Y}_{t}(\omega) y +  A^{Z}_{t}(\omega) z$. In this case, $\Lambdab=\{(A^{Y},A^{Z})\}$ and therefore
$$
\Xc_{0}(l)=E\left[ \int_{0}^{T} L_{s}  \tilde g(s,A^{Y}_s,A^{Z}_s) ds + L_{T} \tilde \Phi(l/L_{T} )\right]\;,
$$  
with $L$ given by 
 $$
L_{t}  =1+\int_{0}^{t}L_{s} A^{Y}_{s}ds+\int_{0}^{t}L_{s} A^{Z}_{s}dW_{s}\;,\;t\in [0,T].
$$ 
Then,  the dual formulation of Proposition \ref{prop dualite existence duale} above drops down to finding $\hat l$ which maximizes $lm-\Xc_{0}(l)$. This generalizes the result of \cite{FoLe99}  and \cite{BoElTo09} obtained  for quantile hedging problems in linear models of financial markets. 
 }
\end{Remark}

\no {\bf Proof of Proposition \ref{prop dualite existence duale}.}  We split the proof in two steps.

{\bf Step  1.} For ease of notations, we set $\hat L:=L^{\hat \lambda}$. By optimality of $\hat l$, one has
$$
 \hat l m- E\left[\hat L_{T}  \tilde \Phi(\hat l/\hat L_{T})\right] \ge  m(\hat l+\iota) - E\left[\hat L_{T} \tilde \Phi((\hat l+\iota)/\hat L_{T} )\right] \;,
$$
for all $\iota>-\hat l$. Since  $\tilde \Phi$ is by construction $\Pas$ convex, this implies that  $\zeta_{\iota} := \nabla \tilde \Phi((\hat l+\iota)/\hat L_{T} )$ satisfies
$
m\iota\le E[{\zeta_{\iota} }] \iota$, for all $\iota>-\hat l$, recall \Hd1 \reff{item: Hd1 C1b} and \reff{item: Hd1 borne nabla}.
Taking $\iota$ of the form $-1/n$ and then $1/n$, for $n\to \infty$, and using  \Hd1 \reff{item: Hd1 C1b} and \reff{item: Hd1 borne nabla} then leads  to 
\be\label{eq: def zeta pour p}
m=   E[\zeta] \;\mbox{ where }\;  \zeta:=\nabla \tilde \Phi(\hat l/\hat L_{T} ).
\ee
We now appeal to \Hd1 \reff{item: Hd1 bi fenchel egale Phi} to deduce that 
\be\label{eq:  zeta pour transfo fenchel p q}
\Phi(\zeta)=\zeta (\hat l/\hat L_{T} )-\tilde \Phi(\hat l/\hat L_{T} ).
\ee
 
By construction, $\tilde \Phi$ is $\Pas$ $1$-Lipschitz and non-decreasing, i.e.~$\zeta \in \Lb_0 ([0,1]) $.
{In view of \reff{eq: def zeta pour p}}, the martingale representation Theorem then implies that we can find $\hat \alpha \in \Ab_{m}$ such that $\hat M_{T}:=M^{m,\hat \alpha}_{T}=\zeta$.  

{\bf  Step 2.} We now write $(\hat \nu,\hat \vc):=\hat \lambda${ and fix $\lambda=(\nu,\theta)\in\Lambdab$ to be chosen later on.}
Clearly,  $\Lambdab$ is convex. Hence,  $\lambda^{\eps}:=(1-\eps)(\hat \nu,\hat \vc)+\eps(\nu,\vc) \in \Lambdab$, $\eps\in [0,1]$. Moreover, direct computations show that  
$$
\frac{\partial}{\partial \eps} L^{\lambda^{\eps}} |_{\eps=0} = \hat L {\hat R}\;\mbox{ where }\;\hat R:= \int_{0}^{\cdot } (\delta \nu_{s} -\delta \vc_{s}\hat \vc_{s})ds + \int_{0}^{\cdot} \delta \vc_{s} dW_{s},
$$
in which we use the notations  $\delta \lambda:=(\delta \nu,\delta \vc):=(\nu-\hat \nu,\vc-\hat \vc)$.

Recalling that elements of $\Lambdab$ take bounded values, see {Remark \ref{rem: domaine tilde Phi et tilde g}}, and  arguing as in Step 1, one easily checks that the optimality condition 
$
{\rm X}_{0}^{\hat l, \lambda^{\eps}}\ge {\rm X}_{0}^{\hat l, \hat \lambda}
$, for all $\eps \in [0,1]$, 
implies that   $\hat \eta:=\nabla \tilde g(\cdot,\hat \lambda)$ satisfies   
\begin{align}
0\le&  E\left[\int_{0}^{T} \hat L_{s}\left(\hat R_{s}\tilde g(s,\hat \lambda_{s})+\hat \eta_{s}^{\top} \delta \lambda_{s}\right)ds 
+\hat R_{T}\hat L_{T}( \tilde \Phi(\hat l/\hat L_{T})\right]\nonumber\\
&+E\left[- (\hat l/\hat L_{T})\nabla \tilde \Phi(\hat l/\hat L_{T}))\right]\nonumber\\
=&E\left[\int_{0}^{T} \hat L_{s}\left(\hat R_{s}\tilde g(s,\hat \lambda_{s})+\hat \eta_{s}^{\top} \delta \lambda_{s}\right)ds 
- \hat R_{T}\hat L_{T} \Phi(\hat M_{T})\right]\;,  \label{eq: calc variation}
\end{align}
{in which we used \reff{eq: def zeta pour p}, \reff{eq:  zeta pour transfo fenchel p q} and the relation $\zeta=\hat M_{T}$ to deduce the second  equality.}
Let $(\hat Y,\hat Z)$ be defined by 
\beq\label{eq: def Y dans verification}
\hat Y:=\hat L^{-1}E_{.}\left[\hat L_{T} \Phi(\hat M_{T})-\int_{.}^{T} \hat L_{s}\tilde g(s,\hat \lambda_{s})ds\right] \;\mbox{ and }\; \hat Z:=\bar Z- \hat Y \hat \vc,
\eeq
 where $\bar Z\in \Hb_{2}$ is implicitly given by 
\beq
\label{eq:dynahetLhatY}
\hat L_{t}\hat Y_{t}=\hat L_{T} \Phi(\hat M_{T})-\int_{t}^{T} \hat L_{s}\tilde g(s,\hat \lambda_{s})ds -\int_{t}^{T} \hat L_{s}\bar Z_{s} dW_{s},\quad 0\le t \le T\;.
\eeq
{The above combined with  \reff{eq: calc variation} implies
\b*
0&\le& E\left[\int_{0}^{T} \hat L_{s}\left(\hat R_{s}\tilde g(s,\hat \lambda_{s})+\hat \eta_{s}^{\top} \delta \lambda_{s}\right)ds 
- \hat R_{T}\hat L_{T} \hat Y_{T}\right].
\e*}
Recalling the definition of $\hat R$ and $\hat \eta$ and {applying It\^{o}'s Lemma, this}  leads to 
\begin{align}
0&\le E\left[\int_{0}^{T} \hat L_{s}\left( \hat \eta_{s}- (\hat Y_{s},\hat Z_{s} )\right)^{\top} \delta \lambda_{s} ds \right]
\nonumber\\
&=
E\left[\int_{0}^{T} \hat L_{s}\left( \nabla \tilde g(s,\hat \lambda_{s})- (\hat Y_{s},\hat Z_{s} )\right)^{\top} \delta \lambda_{s} ds\right].\label{eq: esp ge 0 avant selec mes}
\end{align}
By Assumption \Hd1 \reff{item: Hd1 C1b}, Remark \ref{rem: domaine tilde Phi et tilde g} and \cite[Theorem 18.19, p.~605]{AliprantisBorder}, one can choose  $\bar \lambda \in \Lambdab$ such that 
\b*
&\bar \lambda= {\rm argmin}\left\{f(\cdot,u,v),\;(u,v)\in \dom{\tilde g(\cdot)}\right\}\;\Leb\x \P{\rm -a.e.}
\e*
where 
$$
f:(\omega,s,u,v)\mapsto\left( \nabla \tilde g( \omega,s,\hat \lambda_{s}(\omega) )- (\hat Y_{s}(\omega) ,\hat Z_{s}(\omega))\right)^{\top}  ( u-\hat \nu_{s}(\omega) ,v-\hat \vc_{s}(\omega) ).
$$
Considering now Relation \eqref{eq: esp ge 0 avant selec mes} with $\lambda$ chosen to be equal to $\bar \lambda\Indi{\{f(\cdot,\bar \lambda)<0\}}$, we see that, for $\Leb \x \P$-a.e. $(\omega,t)\in \Omega \x [0,T]$, the gradient   $\Delta_{t}(\omega)$ at   $\hat \lambda_{t}(\omega)$  of the convex map 
$$
(u,v)\in \dom{\tilde g(\omega,t,\cdot)}\mapsto F(\omega,t,u,v):=  \tilde g(\omega,t,u,v)-u\hat Y_{t}(\omega)-v^{\top} \hat Z_{t}(\omega)
$$
satisfies
$$
\Delta_{t}(\omega)^{\top}(b-\hat \lambda_{t}(\omega))\ge 0\;, \quad \;\mbox{ for all } b \in  \dom{\tilde g(\omega,t,\cdot)} .
$$ 
This implies that $\hat \lambda_{t}(\omega)$ minimizes $F(\omega,t,\cdot)$ for $\Leb \x \P$-a.e. $(\omega,t)\in \Omega \x [0,T]$ and therefore we compute
$$
\tilde g(\cdot,\hat \lambda)=\hat \lambda^{\top}(\hat Y ,\hat Z)-g(\cdot,\hat Y,\hat Z)\;\,\,\Leb\x \P-a.e.
$$
by \Hd1 \reff{item: Hd1 bi fenchel egale g}. 
Combining the above identity with \reff{eq:dynahetLhatY} leads to $(\hat Y,\hat Z)=(Y^{m,\hat \alpha},Z^{m,\hat \alpha})$. Then, by using \reff{eq: def zeta pour p}, \reff{eq:  zeta pour transfo fenchel p q} and \reff{eq: def Y dans verification}, in which $\hat L_{0}=1$, we obtain 
\b*
Y^{m,\hat \alpha}_{0}&=&E\left[\hat L_{T} \Phi(\hat M_{T})-\int_{0}^{T} \hat L_{s}\tilde g(s,\hat \lambda_{s})ds\right]\\
&=& E\left[\hat L_{T}\left(\zeta \hat l/\hat L_{T}- \tilde \Phi(\hat l/\hat L_{T})\right)-\int_{0}^{T} \hat L_{s}\tilde g(s,\hat \lambda_{s})ds\right]\\
&=&\hat l m- E\left[\hat L_{T}\tilde \Phi(\hat l/\hat L_{T})+\int_{0}^{T} \hat L_{s}\tilde g(s,\hat \lambda_{s})ds\right].
\e*
In view of Proposition \ref{prop: borne inf dualite}, this concludes the proof.
\ep
 \vs2
 
 We now state the reciprocal statement: existence in the primal problem provides existence in the dual one. Here again, we need to impose some additional technical conditions. 
 \vs2

{\bf Assumption \Hd2}  The following holds for 
$\Leb \x \P$-a.e.  $(t,\omega) \in [0,T]\x \Omega$:
\begin{enumerate}[(a)]
\item\label{item: Hd2 C1b} the maps $   \Phi(\omega,\cdot )$ and $  g(\omega,t,\cdot)$ are $C^{1}_{b}$ on $[0,1]$ and $\R\x \R^{d}$ respectively;
\item\label{item: Hd2 borne nabla}  $|\nabla    \Phi(\omega,\cdot)|\le \chi_{\Phi}(\omega)$, for some $\chi_{\Phi}\in \Lb_{2}(\R)$.
\end{enumerate}
\vs2

\begin{Proposition}\label{prop dualite existence primal}  Let Assumption \Hd2 hold.  Let $l>0$ be fixed and assume   that there exists $\hat m\in [0,1]$ and $ \hat \alpha \in \Ab_{\hat m}$  such that
\be
&\sup\limits_{m\in [0,1]}  \sup\limits_{\alpha \in \Ab_{m}} \left(ml - \Yc_{0}(m)\right)=  \hat m l-  Y^{\hat m,\hat \alpha}_{0}.&
 \ee
 Then,  there exists $\hat \lambda\in \Lambdab$ such that 
 $$
  \Yc_{0} (\hat m)   = \hat m l-    \Xc_{0}(  l)= \hat m l-  {\rm  X}_{0}^{  l,\hat  \lambda} \;,
 $$ 
and $\hat \lambda$ satisfies \reff{eq: prop dualite existence duale - relation variables} with $m=\hat m$ and $\hat l=l$.
  
 \end{Proposition}
 
 \proof Given $\eps\in [0,1]$, a martingale $M$ with values in $[0,1]$, $m:=M_{0}$, we set $m_{\eps}:=\hat m +\eps(m-\hat m)$,  $M^{\eps}:=\hat M +\eps (M-\hat M)$, where $\hat M:=M^{\hat m,\hat \alpha}$. For ease of notation, we set  $(\hat Y,\hat Z):=(Y^{\hat m,\hat \alpha},Z^{\hat m,\hat \alpha})$ and denote by $(Y^{\eps},Z^{\eps})$ the solution of $\BSDE(g,\Phi(M^{\eps}_{T}))$,  $\delta m:=m-\hat m$, $(\delta M,\delta Y^{\eps}, \delta Z^{\eps}):=(M-\hat M,Y^{\eps}-\hat Y,Z^{\eps}-\hat Z)$.
 
 \no{\bf Step 1. } {\sl We first show that  $\eps^{-1}(\delta Y^{\eps}_{s},\delta Z^{\eps}_{s})$ converges in $\Sb_{2}\x \Hb_{2}$ as $\eps\to 0$ to the solution $(\nabla Y,\nabla Z)$ of  
  \beq\label{eq: def bsde tangente en eps}
 \nabla Y_{t}= \nabla \Phi(\hat M_{T})\delta M_{T}+ \int_{t}^{T}\nabla g(s,\hat Y_{s},\hat Z_{s})^{\top}(\nabla Y_{s},\nabla Z_{s})ds -\int_{t}^{T} \nabla Z_{s} dW_{s} .
 \eeq
 }
 
 First note that existence and uniqueness of the solution to the above BSDE in guaranteed by Assumption \Hd2. 
 
Letting $\xi^\eps:=\eps^{-1}(\Phi(M_T^\eps)-\Phi(\hat M_T))$, one easily checks that  $\eps^{-1}(\delta Y^{\eps}_{s},\delta Z^{\eps}_{s})$  solves 
 $$ \frac{\delta Y^{\eps}_{s}}{\eps} = \xi^\eps - \int_s^T \frac{\delta Z^{\eps}_{r}}{\eps} dW_r + \int_s^T \left(A_{r}^{Y,\eps} \frac{\delta Y^{\eps}_{r}}{\eps} + A_{r}^{Z,\eps} \frac{\delta Z^{\eps}_{r}}{\eps}\right) dr,$$
where 
$$
 A_{r}^{Y,\eps}:= \int_0^1 \partial_y g(r,\hat{Y}_r +\theta \delta Y^{\eps}_r, \hat{Z}_r) d\theta
 \;\mbox{ and }\;
 A_{r}^{Z,\eps}:= \int_0^1 \partial_z g(r,Y^\eps_r, \hat{Z}_r+\theta \delta Z^{\eps}_r) d\theta.
 $$
In the above, $ \partial_y g$ and $ \partial_z g$ denotes respectively the partial gradients of $g$ with respect to $y$ and $z$, recall \Hd2. 
The Assumption \Hg~ implies $|A^{Y,\eps}|+|A^{Z,\eps}|\leq K_{g}$.  

   We now  set $U^{\eps}:=\eps^{-1} \delta Y^{\eps}_{s}-\nabla Y$, $V^{\eps}:=\eps^{-1} \delta Z^{\eps}_{s}-\nabla Z$ and $\zeta^{\eps}:=\xi^\eps-\nabla \Phi(\hat{M}_T) \delta M$. The pair $(U^{\eps},V^{\eps})$ is an element of $\Sb_{2}\x \Hb_{2}$ and solves 
$$ 
U^{\eps}_{s} = \zeta^{\eps}- \int_s^T V^{\eps}_r dW_r + \int_s^T \left(A_{r}^{Y,\eps} U_r^\eps + A_{r}^{Z,\eps} V^{\eps}_r +R^{\eps}_r\right) dr\;,\quad 0\le s\le T \;,
$$
with
$$ 
R^{\eps}_r:= \nabla Z_r (A_{r}^{Z,\eps}-\partial_z g(r,\hat{Y}_r,\hat{Z}_r))+\nabla Y_r (A_{r}^{Y,\eps}-\partial_y g(r,\hat{Y}_r,\hat{Z}_r)) \;,\quad 0\le r\le T \;.
$$
Hence, by stability for Lipschitz BSDEs (see Proposition \ref{PropStability} in the Appendix) there exists a constant $C>0$ (which does not depend on $\eps$) such that
\begin{equation}
\label{eq:diff1}
\|U^{\eps}\|_{\Sb_{2}}^2 + \|V^{\eps}\|_{\Hb_{2}}^2 \leq C \left( \|\zeta^{\eps}\|_{\Lb_{2}}^2 +\|R^{\eps}_r\|^2_{\Hb_{2}} \right).
\end{equation}
The result of Step 1. will follow if we prove that the right-hand side of the inequality \eqref{eq:diff1} vanishes as $\eps$ tends to zero. 
The convergence of $\|R^{\eps}_r\|^2_{\Hb_{2}}$ to $0$ follows from Assumption \Hd2 and the convergence of $M_T^\eps$ to $M_T$. As for the second term, it suffices to prove   that $(Y^\eps,Z^\eps)_\eps$ converges in $\Sb_{2}\x \Hb_{2}$ to $(\hat{Y},\hat{Z})$, and to appeal to \Hg~and \Hd2.  The latter is obtained by standard stability results,  see Proposition \ref{PropStability} below, which imply the existence of a constant $C>0$ (which does not depend on $\eps$) such that
\begin{align*}
\|Y^\eps-\hat{Y}\|_{\Sb_{2}}^2 + \|Z^\eps-\hat{Z}\|_{\Hb_{2}}^2 &\leq C \|\Phi(M_T^\eps)-\Phi(\hat{M}_T)\|_{\Lb_{2}}^2\  \longrightarrow_{\eps \to 0} 0 .
\end{align*} 
In the latter,  the convergence follows from Lebesgue's dominated convergence Theorem and assumption \Hd2.
 
 \no {\bf Step 2. } By optimality of  $(\hat m,\hat \alpha)$, $ Y^{\eps}_{0}-m_{\eps}l- {\hat Y_{0}}+\hat m   l\ge 0$, for any $\eps>0$.
 In view of Step 1, dividing by $\eps>0$ and sending $\eps\to 0$ leads to 
\begin{align*}
0&\le \nabla \Phi(\hat M_{T})\delta M_{T}  -l \delta m  +\int_{0}^{T}\nabla g(s,\hat Y_{s},\hat Z_{s})^{\top}(\nabla Y_{s},\nabla Z_{s})ds - \int_{0}^{T} \nabla Z_{s} dW_{s}\\
& =\nabla Y_{0} -l \delta m ,
\end{align*}
after possibly passing to a subsequence.

Set $\hat L:=L^{\hat \lambda}$ where  $\hat\lambda:=\nabla g (\cdot,\hat Y ,\hat Z)$. Observe that the latter  belongs to $\Lambdab$. For later use, also notice that  
\be\label{eq: hat beta comme gradient preuve 2 dualite}
g(\cdot,\hat Y,\hat Z)= (\hat \nu,\hat \vc)^{\top}(\hat Y,\hat Z)-\tilde g (\cdot,\hat \nu,\hat \vc), 
\ee
see e.g. \cite{Rockafellar}. Then, it follows from \reff{eq: def bsde tangente en eps} that $\hat L \nabla Y$ is a martingale. The previous inequality thus implies that 
$$
0\le \hat L_{0} \nabla Y_{0}-l\delta m= E\left[\hat L_{T} \nabla Y_{T}\right]-l\delta m= E\left[\hat L_{T} \delta M_{T}\left(\nabla \Phi(\hat M_{T}) -l/\hat L_{T}  \right) \right],
$$
in which we used the fact that $\hat L_{0}=1$ and $E[\delta M_{T}]=\delta m$. Since $M_{T}$ can be any arbitrary random variable with values in $[0,1]$, this shows that, $\Pas$, $\hat M_{T}(\omega)$ minimizes $m\in [0,1]\mapsto \Phi(\omega,m)-ml/\hat L_{T}(\omega)$.  Hence, 
$$
\hat M_{T} l   - \hat L_{T}\Phi(\hat M_{T})=\hat L_{T}\tilde \Phi(l/\hat L_{T}),
$$
see e.g. \cite{Rockafellar}. Combining the above identity together with \reff{eq: hat beta comme gradient preuve 2 dualite} and using It\^{o}'s Lemma leads to 
$
l\hat m- \hat Y_{0}={\rm X}_{0}^{\hat l,\hat \lambda}.  
$ 
One concludes by appealing to Proposition  \ref{prop: borne inf dualite}.
\ep
  
\section{Proof of Theorem \ref{thm : BSDE characterization bar Yc} }\label{SEC: proof thm BSDE characterization bar Yc}

{In all this section, we use the notations introduced at the beginning of Section \ref{SEC: main characterization}.} The first main result provides a dynamic programming principle for the family $\{\Yc_{\tau}^{\alpha},\;\tau \in \Tc, \alpha \in \Ab_{0}\}$.

\begin{Proposition}\label{prop : DPP} For all $(\tau_{1},\tau_{2},\alpha)\in \Tc\x \Tc\x \Ab_{0}$ such that  $\tau_{1}\le \tau_{2}$, we have 
\b*
\Yc_{\tau_{1}}^{\alpha} &=& \essinf_{\alpha'\in \Ab_{\tau_{1}}^{\alpha} } \Eg{\tau_{1},\tau_{2}} {\Yc_{\tau_{2}}^{\alpha'}} \;.
\e*
\end{Proposition}

\proof  We prove the two corresponding  inequalities separately. \\
\no{\bf Step 1. $\Yc_{\tau_{1}}^{\alpha}\ge \essinf_{\alpha'\in \Ab_{\tau_{1}}^{\alpha} } \Ec^g_{\tau_{1},\tau_{2}}\left[ \Yc_{\tau_{2}}^{\alpha'}\right]$}.\\
It follows from Lemma \ref{lemma repres Y} below that  there exists $(\alpha^n)_n$ in $\Ab_{\tau_{1}}^{\alpha}$ such that the sequence $(\Ec^g_{\tau_1,T}[\Phi(M_T^{\alpha^n})])_n$ is non-increasing and
\be\label{eq conv nnn} \lim_{n\to \infty} \Ec^g_{\tau_1,T}[\Phi(M_T^{\alpha^n})] &=& \Yc_{\tau_{1}}^{\alpha}, \quad \Pas \ee
Since    $\alpha^n \in \Ab_{\tau_{2}}^{\alpha^n}$  for every $n\geq 1$, we deduce that
\b* \Yc_{\tau_{2}}^{\alpha^n} &\leq& \Ec^g_{\tau_2,T}[\Phi(M_T^{\alpha^n})]\,.\e*
By comparison for BSDEs with Lipschitz continuous drivers on the time interval $[\tau_1,\tau_2]$, this implies
\b* \Ec^g_{\tau_1,\tau_2}[\Yc_{\tau_{2}}^{\alpha^n}]  &\leq&  \Ec^g_{\tau_1,\tau_2}[\Ec^g_{\tau_2,T}[\Phi(M_T^{\alpha^n})]]=\Ec^g_{\tau_1,T}[\Phi(M_T^{\alpha^n})] \,,\e*
leading to 
\b* \essinf_{\alpha'\in \Ab_{\tau_{1}}^{\alpha} } \Ec^g_{\tau_{1},\tau_{2}}\left[ \Yc_{\tau_{2}}^{\alpha'}\right] & \leq & \Ec^g_{\tau_1,T}[\Phi(M_T^{\alpha^n})] \,,\e*
Letting $n$ go to infinity in the above inequality, \reff{eq conv nnn}  provides directly
\b* \essinf_{\alpha'\in \Ab_{\tau_{1}}^{\alpha} } \Ec^g_{\tau_{1},\tau_{2}}\left[ \Yc_{\tau_{2}}^{\alpha'}\right] &\leq& \Yc_{\tau_{1}}^{\alpha} \,.\e*

\no{\bf Step 2. $\Yc_{\tau_{1}}^{\alpha}\le \essinf_{\alpha'\in \Ab_{\tau_{1}}^{\alpha} } \Ec^g_{\tau_{1},\tau_{2}}\left[ \Yc_{\tau_{2}}^{\alpha'}\right]$}.\\

Fix $\alpha'$ in $\Ab_{\tau_{1}}^{\alpha}$.
Lemma \ref{lemma repres Y} below ensures the existence of a sequence $(\alpha'_n)_n$ in $\Ab_{\tau_{2}}^{\alpha'}$ such that $(\Ec^g_{\tau_2,T}[\Phi(M_T^{\alpha'^n})])_n$ is non-increasing and
$$ \lim_{n\to \infty} \Ec^g_{\tau_2,T}[\Phi(M_T^{\alpha'_n})]= \Yc_{\tau_{2}}^{\alpha'}, \quad \Pas$$
In view of Remark \ref{rem: borne unif Yc}, the convergence holds in $\Lb_{2}$ as well.  Thus the stability result of Proposition \ref{PropStability} below indicates that 
$\Ec^g_{\tau_1,T}[\Phi(M_T^{\alpha'_n})]$ converges to $\Ec^g_{\tau_1,\tau_2}[\Yc_{\tau_{2}}^{\alpha'}]$ in $\Lb_2$. 
In addition,   $\alpha'_n \in \Ab_{\tau_{2}}^{\alpha'} \subset \Ab_{\tau_{1}}^{\alpha}$ by construction. Combining the above leads to 
\b* \Ec^g_{\tau_1,\tau_2}[\Yc_{\tau_{2}}^{\alpha'}] &=& \lim_{n\to \infty} \Ec^g_{\tau_1,T}[\Phi(M_T^{\alpha'_{n}})] \;\geq\; \Yc_{\tau_{1}}^{\alpha} \; . \e*
The arbitrariness of $\alpha'\in\Ab_{\tau_{1}}^{\alpha}$ allows one to conclude
\b* \essinf_{\alpha'\in \Ab_{\tau_{1}}^{\alpha} } \Ec^g_{\tau_1,\tau_2}[\Yc_{\tau_{2}}^{\alpha'}] &\geq& \Yc_{\tau_{1}}^{\alpha} \;. \e*
\ep

\begin{Lemma}\label{lemma repres Y} Fix $\theta,\tau \in \Tc$, {with $\theta\ge \tau$,} $\mu\in \Lb_{0}([0,1],\Fc_{\tau})$ and $\alpha \in \Ab_{\tau,\mu}$. Then, there exists a sequence $(\alpha'_{n})\subset \Ab_{\tau,\mu}^{\theta, \alpha}$ :=$\{\alpha' \in \Ab_{\tau,\mu}, \; \alpha'\textbf{1}_{[0,\theta)}=\alpha \textbf{1}_{[0,\theta)} \}$such that $\lim_{n} \downarrow \Eg{\theta,T}{\Phi(M^{\tau,\mu,\alpha_{n}'}_{T})}={\Yc^{\alpha}_{\theta}(M^{\tau,\mu,\alpha}_{\theta})}$ $\Pas$ 
\end{Lemma}
\proof It suffices to show that the family  $\{J(\alpha'):=\Eg{\theta,T}{\Phi(M^{\tau,\mu,\alpha'}_{T})},\;\alpha'\in \Ab_{\tau,\mu}^{\theta,\alpha}\}$ is directed downward, see e.g. \cite{Ne75}. 
Fix $\alpha'_{1}, \alpha'_{2}$ in $\Ab_{\tau,\mu}^{\theta,\alpha}$ and set 
$$
\tilde{\alpha}':=\alpha\1_{[0,\theta)} +\1_{[\theta,T]} (\alpha'_{1} \1_{A}+\alpha'_{2} \1_{A^{c}} )
$$
where $A:=\{J(\alpha'_{1})\le J(\alpha'_{2})\}\in \Fc_{\theta}$, so that $\tilde{\alpha}'\in \Ab_{\tau,\mu}^{\theta,\alpha}$ and 
\b*
J(\tilde \alpha')=\Eg{\theta,T}{\Phi(M^{\tau,\mu,\alpha'_{1}}_{T})\1_{A}+\Phi(M^{\tau,\mu,\alpha'_{{2}}}_{T})\1_{A^{c}}}=\min\{J(\alpha'_{1}{)},J(\alpha'_{2}{)}\}.
\e*
\ep
   
We now observe that the family $(\Yc^{\alpha})_{\alpha\in\Hb_{2}}$ is  {l\`adl\`ag on countable sets}.
If in addition $\Phi$ is assumed to be continuous, the process $(\Yc^{\alpha})_{\alpha\in\Hb_{2}}$ is even indistinguishable from a c\`adl\`ag process.  

\begin{Proposition}\label{prop : aggregation cadlag} Fix $\alpha \in \Ab_{0}$. Then, $\Yc^\alpha$ is {l\`adl\`ag  on countable sets}. Besides, if $m\in [0,1]\mapsto \Phi(\omega,m)$ is continuous for  $\P$-a.e. $\omega \in \Omega$, then $\Yc^\alpha$ is indistinguishable from a c\`adl\`ag process.
\end{Proposition}

\proof
Fix $\alpha\in\Ab_{0}$. Proposition \ref{prop : DPP} and Remark \ref{rem: borne unif Yc} imply that $-\Yc^{\alpha}$ is a $-g(-\cdot)$-supermartingale in the sense of \cite{ChPe00} (a $g$-submartingale in the sense of \cite{Pe99}). It follows from the non-linear up-crossing Lemma, see \cite[Theorem 6]{ChPe00}\footnote{Note that   \cite[Theorem 6]{ChPe00} restricts to  positive $g$-supermartingales. However, the proof can be reproduced without difficulty under the integrability condition of Remark \ref{rem: borne unif Yc}. In addition, \cite[Theorem 6]{ChPe00} implies that $E^{\mathbb{Q}}[D_a^b(\Yc^{\alpha},n)] \leq \Yc^{\alpha}_0\wedge b \leq b$, where $D_a^b(\Yc^{\alpha},n)$ denotes the number of down crossing of $\Yc^{\alpha}$ from an interval $[a,b]$ on a discrete time-grid  $0=t_0\leq t_1 \leq \cdots \leq t_n=T$ and $\mathbb{Q}$ is a particular measure absolutely continuous with respect to $\P$. To conclude, it is enough to reproduce the proof of \cite[Chapter VI Theorem (2) point 1)]{DellacherieMeyer}.}, 
that the following limits
$$ \lim_{s\in {{\rm D}}\cap  (t,T]\downarrow t} \Yc^{\alpha}_{s} \textrm{ and } \lim_{s\in {{\rm D}}\cap  [0,t)\uparrow t} \Yc^{\alpha}_{s} $$ 
are well-defined for every $t$ in $[0,T]$, $\Pas$, {and for all countable set ${\rm D}$. So} is the process 
$${\bar \Yc}^{\alpha}_{t}:=\lim_{s\in \Q\cap  (t,T]\downarrow t} \Yc^{\alpha}_{s}\,, \qquad t\in [0,T]\;.$$

Besides, ${\bar \Yc}^{\alpha}$ is by definition c\`ad. Assuming that $\Phi$ is continuous, we will prove that, for every stopping time $\tau$, it holds that: 

\begin{equation}
\label{eq:cadlagtemp1}
{\bar \Yc}^{\alpha}_\tau = \essinf_{\alpha'\in \Ab_{\tau}^{\alpha} } \Ec^g_{\tau,T}\left[ \Phi(M^{\alpha'}_T)\right](= \Yc^{\alpha}_\tau) \qquad \Pas
\end{equation}
By \cite[Chapter IV. (86), p.~220]{DellacherieMeyer},  
{the r}elation \eqref{eq:cadlagtemp1} entails that $\Yc^\alpha$ and $\bar \Yc^\alpha$ are undistinguishable showing that $\Yc^{\alpha}$ is undistinguishable from a c\`adl\`ag process. The rest of the proof is devoted to prove   \eqref{eq:cadlagtemp1}. 
 
For this purpose, let us introduce $(\tau_{n})_{n}$, a decreasing sequence of stopping times with values in $[0,T]\cap \Q$  such that $\tau\le \tau_{n}\le \tau+ n^{-1}$ and ${\bar \Yc}_{\tau}^{\alpha}=\lim_{n\to \infty} \Yc_{\tau_{n}}^{\alpha}$. \\\\
\noindent
\no{\bf Step 1. ${\bar \Yc}_{\tau}^{\alpha}\le \essinf_{\alpha'\in \Ab_{\tau}^{\alpha} } \Ec^g_{\tau,T}\left[ \Phi(M^{\alpha'}_T)\right]$}.\\
{\bf a.} Fix $\alpha'\in \Ab_{\tau}^{\alpha}$ and set  
$$ \lambda_n:=\left( \frac{M_{\tau_{n}}^\alpha}{M_{\tau_{n}}^{\alpha'}} \wedge \frac{1-M_{\tau_{n}}^\alpha}{1-M_{\tau_{n}}^{\alpha'}} \right) \textbf{1}_{\{M_{\tau_{n}}^\alpha\notin \{0,1\}\}} \in [0,1],$$
with the convention $a/0=\infty$ for $a>0$. 
Using the fact that $M^{\alpha'}_{\tau_{n}}+\int_{\tau_{n}}^{T }\alpha'_{s}dW_{s}=M^{\alpha'}_{T}\in [0,1]$, direct computations lead to
$$
 0\le M^{\alpha}_{\tau_{n}} -\lambda_{n} M^{\alpha'}_{\tau_{n}}
\le 
M^{\alpha}_{\tau_{n}}+\lambda_{n}\int_{\tau_{n}}^{T}\alpha'_{s}dW_{s}
\le 
 M^{\alpha}_{\tau_{n}} +\lambda_{n}(1- M^{\alpha'}_{\tau_{n}})\le 1.
$$
We set $\alpha^{'}_{n}:=\alpha \1_{[0,\tau_{n})}+ \lambda_n \alpha' \1_{[\tau_{n},T]}$. The above implies that  $\alpha'_{n}$ belongs to $\Ab_{\tau_{n}}^{\alpha}$. 

\noindent {\bf b.} Now we prove that $M_T^{\alpha'_n}$ converges $M^{\alpha'}_T$ in $\Lb_2$ as $n$ goes to infinity, possibly up to a subsequence. Since both have norms bounded by $1$, it suffices to show the $\Pas$ convergence, possibly up to a subsequence.  To see this, first note that 
$$M_T^{\alpha'_n}-M^{\alpha'}_T = M_{\tau_{n}}^\alpha - M_{\tau_{n}}^{\alpha'} + \int_{\tau_{n}}^T (\lambda_n-1) \alpha'_s dW_s,$$
from which we deduce that
\begin{align*}
M_T^{\alpha'_n}-M^{\alpha'}_T  &= M_{\tau_{n}}^\alpha - M_{\tau_{n}}^{\alpha'} + (\lambda_n-\textbf{1}_{\{M_{\tau_{n}}^\alpha\notin \{0,1\}\}})\; \int_{\tau_{n}}^T \alpha'_s dW_s\\
& - \textbf{1}_{\{M_{\tau_{n}}^\alpha \in \{0,1\}\}} \; \int_{\tau_{n}}^T \alpha'_s dW_s .
\end{align*}
Since $\tau_{n}\to \tau$ $\Pas$ and $\alpha'=\alpha$ on $[\![0,\tau]\!]$, the above  construction implies that $\lim_{n\to \infty} M_{\tau_{n}}^\alpha - M_{\tau_{n}}^{\alpha'}=0  \; \Pas$ and $\lim_{n\to \infty} \lambda_n =\lim_{n\to \infty} \textbf{1}_{\{M_{\tau_{n}}^\alpha\notin \{0,1\}\}}  \; \Pas$ 
It thus only remains to prove that $\textbf{1}_{\{M_{\tau_{n}}^\alpha \in \{0,1\}\}} \; \int_{\tau_{n}}^T \alpha'_s dW_s  \to 0$ $\Pas$ First note  that $\alpha'\1_{[\tau_{n},T]}=0$ on $\{M^{\alpha'}_{\tau_{n}} \in \{0,1\}\}$. This follows from the martingale property of this process with values in $[0,1]$. Hence, it suffices to consider $\textbf{1}_{\{M_{\tau_{n}}^{\alpha'} \ne M_{\tau_{n}}^\alpha \in \{0,1\}\}} \; \int_{\tau_{n}}^T \alpha'_s dW_s$.
But, since $M_{\tau}^{\alpha'} = M_{\tau}^\alpha $, 
\begin{align*}
\P[M_{\tau_{n}}^{\alpha'} \ne M_{\tau_{n}}^\alpha \in \{0,1\}] 
&\leq 
\P[M_{\tau_{n}}^{\alpha'} \ne M_{\tau_{n}}^\alpha  ] 
= 
\P\left[\left|\int_\tau^{\tau_{n}} (\alpha_s-\alpha_s') dW_s\right|>0\right]\\
& \rightarrow_{n\to \infty} 0.
\end{align*}

\noindent {\bf c.} Now, since $\Phi$ is continuous {and $M^{\alpha'_n}_{T}\in \Lb_{0}([0,1])$},  we get that $\Phi(M^{\alpha'_n}_{T})\to \Phi(M^{\alpha'}_{T})$ in $\Lb_2${, after possibly passing to a subsequence}. The stability property for Lipschitz BSDEs given in Proposition \ref{PropStability} implies that 
\begin{equation}
\label{eq:tempconv1}
\left\| \Ec^g_{\tau^n,T}\left[ \Phi(M^{\alpha'_n}_{T})\right] - \Ec^g_{\tau^n,T}\left[\Phi( M^{\alpha'}_{T})\right]\right\|_{\Lb_2} \rightarrow_{n\rightarrow\infty} 0 \;.
\end{equation} 
 {On the other hand, the bound of Remark \ref{rem: borne unif Yc} implies that}
 \begin{equation}
\label{eq:tempconv2}
\left\| \Ec^g_{\tau^n,T}\left[\Phi( M^{\alpha'}_{T})\right]-\Ec^g_{\tau,T}\left[\Phi( M^{\alpha'}_{T})\right] \right\|_{\Lb_{2}} \rightarrow_{n\rightarrow\infty} 0 \;,
\end{equation} 
by Lebesgue's dominated convergence Theorem and by continuity of the process $\Ec^g_{\cdot,T}\left[\Phi( M^{\alpha'}_{T})\right]$. 
Combining \eqref{eq:tempconv1} and \eqref{eq:tempconv2} leads to 
$$
{\bar \Yc}_{\tau}^{\alpha}=\lim_{n\to \infty} \Yc_{\tau^{n}}^{\alpha}\le \lim_{n\to \infty} \Ec^g_{\tau^{n},T}\left[ \Phi(M^{\alpha'_n}_{T})\right]=\Ec^g_{\tau,T}\left[\Phi( M^{\alpha'}_{T})\right]. 
$$
We conclude by arbitrariness of $\alpha'\in \Ab_{\tau}^{\alpha}$.

\no{\bf Step 2. ${\bar \Yc}_{\tau}^{\alpha}\ge \essinf_{\alpha'\in \Ab_{\tau}^{\alpha} } \Ec^g_{\tau,T}\left[ \Phi(M^{\alpha'}_T)\right]$}.\\
Applying on $[\tau,\tau^n]$ the stability result of Proposition \ref{PropStability}  for the BSDEs with parameters $({\bar \Yc}^\alpha_{\tau},0)$ and $(\Yc^\alpha_{\tau^n},g\1_{[0,\tau^{n})})$, we get
\b*
\left\| {\bar \Yc}^\alpha_{\tau}  - \Ec^g_{\tau,\tau^n}\left[\Yc^\alpha_{\tau^n}\right]\right\|_{\Lb_2} 
& \le &  
C \left(\left\| {\bar \Yc}^\alpha_{\tau}  -  \Yc^\alpha_{\tau^n} \right\|_{\Lb_2} 
+
 E\left[ \int_{\tau}^{\tau^n} \left|g(s,{\bar \Yc}^\alpha_{\tau},0) \right|^2 ds \right] \right)\\
& \le &  
C \left\| {\bar \Yc}^\alpha_{\tau}  -  \Yc^\alpha_{\tau^n} \right\|_{\Lb_2} 
+
 \frac{C}{n} \;, \qquad n\in\N\;, 
\e*
 for some $C>0$, since {the bound of Remark \ref{rem: borne unif Yc} holds for } ${\bar \Yc}^\alpha_{\tau}$, recall that Assumption \Hg~ is in force. Therefore, $\Ec^g_{\tau,\tau^n}[ \Yc^\alpha_{\tau^n}]$ converges to ${\bar \Yc}^\alpha_{\tau}$ as $n$ goes to infinity. Proposition \ref{prop : DPP}  {implies
$
\Ec^g_{\tau,\tau^n}\left[ \Yc^\alpha_{\tau^n}\right] \ge \Yc^\alpha_{\tau}
$. Passing to the limit   leads to the required inequality:}  ${\bar \Yc}_{\tau}^{\alpha}$ $\ge$ $ \Yc^\alpha_{\tau}$ $=$ $ \essinf_{\alpha'\in \Ab_{\tau}^{\alpha} } \Ec^g_{\tau,T}\left[ \Phi(M^{\alpha'}_T)\right]$.
\ep
 \\

In the rest of this section, we complete the proof of Theorem \ref{thm : BSDE characterization bar Yc}. \\

\noindent \textbf{{Proof of Theorem \ref{thm : BSDE characterization bar Yc}.}} Items (i) and (ii) are already proved in Proposition \ref{prop : DPP} and Proposition \ref{prop : aggregation cadlag}, it remains to prove (iii) and (iv).
For $\alpha\in \Ab_{0}$, it follows from Proposition \ref{prop : DPP},  Proposition \ref{prop : aggregation cadlag}  and standard comparison results for BSDEs that $\Yc^\alpha$ is a c\`adl\`ag strong $g$-submartingale in the sense of \cite{Pe99}.  Hence, the existence of a process $( \Zc^{\alpha},\Kc^{\alpha}) \in\Hb_{2}\x \Kb_{2}$ such that \reff{eq : BSDE de la cond minimale - BSDE}  holds   follows from \cite[Theorem 3.3]{Pe99}. We now verify successively that the family $(\Yc^{\alpha},\Zc^{\alpha}, \Kc^{\alpha})_{\alpha \in \Hb_{2}}$ satisfies {\reff{eq : BSDE de la cond minimale - borne uniforme},  \reff{eq : BSDE de la cond minimale - cond de minimalite}, \reff{eq : BSDE de la cond minimale - indep alpha futur}} and the uniqueness of solution for \reff{eq : BSDE de la cond minimale - borne uniforme}-\reff{eq : BSDE de la cond minimale - BSDE}-\reff{eq : BSDE de la cond minimale - cond de minimalite}-\reff{eq : BSDE de la cond minimale - indep alpha futur}.\\

{The bound \reff{eq : BSDE de la cond minimale - borne uniforme} follows directly from  Remark \ref{rem: borne unif Yc} and the representation Theorem 3.3 in \cite{Pe99}, note that   the driver function $g$ does not depend on $\alpha\in\Ab_{0}$.}\\

\no {\bf Step 1. The irrelevance of future property \reff{eq : BSDE de la cond minimale - indep alpha futur}}\\
 For $(\alpha,\tau)\in\Ab_{0}\times\Tc$, observe that $\Ab_{.}^{\alpha'}=\Ab_{.}^{\alpha}$ on $[0,\tau]$ when $\alpha' \in \Ab^{\alpha}_{\tau}$. The definition of $\Yc$ thus implies 
that $
\Yc^{\alpha}\1_{[0,\tau]} $ $=$ $ \Yc^{\alpha'}\1_{[0,\tau]}$ for $ \alpha' \in \Ab^{\alpha}_{\tau}$.
 Hence \reff{eq : BSDE de la cond minimale - indep alpha futur} follows from the uniqueness of the representation provided in    \cite[Theorem 3.3]{Pe99}.\\

\no {\bf Step 2. The minimality property \reff{eq : BSDE de la cond minimale - cond de minimalite}}\\
We follow the arguments in the proof   \cite[Theorem 4.6]{sontouzha11}.
 We fix $(\alpha,\tau_1,\tau_2)\in\Hb_{2}\times\Tc\times\Tc$ such that $\tau_1 \le \tau_2$. For any $\alpha' \in \Ab^{\alpha}_{\tau_1}$, we denote by $(Y^{\alpha'},Z^{\alpha'})$ the solution of the classical BSDE
 \b*
 Y^{\alpha'}_t &=& \Phi(M^{\alpha'}_T) + \int_{t}^T g(s,Y^{\alpha'}_s,Z^{\alpha'}_s)ds -  \int_{t}^T Z^{\alpha'}_s dW_s \;, \qquad 0\le t \le T\;.
 \e*
  Let $L^{\alpha'}$ be the process whose dynamics is given by
  \b*
  L^{\alpha'}_t &=& \exp\left(  \int_{\tau_1}^t \Lambda^z_s dW_s {+} \int_{\tau_1}^t \left(\Lambda^y_s {-} \frac{|\Lambda^z_s|^2}{2}\right) ds  \right) \;, \qquad {\tau_1}\le t \le T \;,
  \e*
  where $(\Lambda^y,\Lambda^z)$ is the linearization process  given by
 \b*
 \Lambda^y &:=& \frac{g(\Yc^{\alpha'}_s,\Zc^{\alpha'}_s)-g(Y^{\alpha'}_s,  {\Zc^{\alpha'}_s})}{\Yc^{\alpha'}_s-Y^{\alpha'}_s}\1_{\{\Yc^{\alpha'}\neq Y^{\alpha'}\}} \;, \\
 \Lambda^z &:=& \frac{g( {Y^{\alpha'}_s}, \Zc^{\alpha'}_s)-g(Y^{\alpha'}_s,Z^{\alpha'}_s)}{|\Zc^{\alpha'}_s-Z^{\alpha'}_s|^2}(\Zc^{\alpha'} - Z^{\alpha'}) \1_{\{\Zc^{\alpha'}\neq Z^{\alpha'}\}} \;. 
  \e*
 This linearization procedure implies that $Y^{\alpha'}_{\tau_1}-\Yc^{\alpha'}_{\tau_1}$ rewrites as
 \be
  Y^{\alpha'}_{\tau_1}- \Yc^{\alpha'}_{\tau_1} &=& 
{E_{\tau_1}\left[ L_{\tau_2}^{\alpha'} (Y^{\alpha'}_{\tau_2}- \Yc^{\alpha'}_{\tau_2})  \right] +} E_{\tau_1}\left[ \int_{\tau_1}^{\tau_2} L^{\alpha'}_s d\Kc^{\alpha'}_s  \right] \nonumber\\
 &\ge& 
 E_{\tau_1}\left[ (\Kc^{\alpha'}_{\tau_2} - \Kc^{\alpha'}_{\tau_1}) \inf_{[\tau_1,\tau_2]} L^{\alpha'}  \right]\label{eqtemp1212}
  \;,
 \ee
 where we used the fact that $Y^{\alpha}-\Yc^{\alpha}\ge 0$. 
 Using H\"older inequality, this implies
 \begin{align*}
& E_{\tau_1}\left[ (\Kc^{\alpha'}_{\tau_2} - \Kc^{\alpha'}_{\tau_1}) \right]^{3}
 \\&\le 
  E_{\tau_1}\left[ (\Kc^{\alpha'}_{\tau_2} - \Kc^{\alpha'}_{\tau_1}) \inf_{[\tau_1,\tau_2]} L^{\alpha'} \right]
  E_{\tau_1}\left[ \sup_{[\tau_1,\tau_2]} (1/L^{\alpha'})\right]
  E_{\tau_1}\left[ (\Kc^{\alpha'}_{\tau_2} - \Kc^{\alpha'}_{\tau_1})^2\right]
 \\
 &\le
 C \;  
 E_{\tau_1}\left[ (\Kc^{\alpha'}_{\tau_2} - \Kc^{\alpha'}_{\tau_1})^2\right]   (Y^{\alpha'}_{\tau_1}-\Yc^{\alpha'}_{\tau_1}) \;,
 \end{align*}
 for some $C>0$ that depends on the uniform bounds on $(\Lambda^y,\Lambda^z)$, recall \Hg. Hence, the estimate \reff{eq : BSDE de la cond minimale - borne uniforme} together with the monotonicity of $\Kc$ implies  
 \be\label{eqtamp111}
 \quad 0\; \le \; E_{\tau_1}\left[ (\Kc^{\alpha'}_{\tau_2} - \Kc^{\alpha'}_{\tau_1}) \right]
 &\le&
 C \eta'_{{\tau_1}} \;\;\;  
 (Y^{\alpha'}_{\tau_1}- \Yc^{\alpha'}_{\tau_1})^{1/3}\;, \qquad \alpha' \in \Ab^{\alpha}_{\tau_1}\;,
 \ee
 where 
 \b*
  \eta'_{{\tau_1}}
  &:=&
  \esssup_{\bar \alpha\in\Ab^{\alpha}_{\tau_1}}E_{\tau_1}\left[ (\Kc^{\bar \alpha}_{\tau_2} - \Kc^{\bar \alpha}_{\tau_1})^2\right]^{1/3} .  
  \e* 
  By the same arguments as in Lemma \ref{lemma repres Y}, we  can find a sequence $(\alpha'_n)_n\subset\Ab^{\alpha}_{\tau_1}$ such that 
  \b*
  \eta'_{{\tau_1}}
  &=&
  \lim_{n\rightarrow\infty} \uparrow  E_{\tau_1}\left[ (\Kc^{\alpha'_n}_{\tau_2} - \Kc^{\alpha'_n}_{\tau_1})^2\right]^{1/3} \;. 
  \e*
 The monotone convergence Theorem together with Jensen's inequality and Relation \reff{eq : BSDE de la cond minimale - borne uniforme} imply that
  \b*
  E [\eta'_{{\tau_1}}]
  &=&
  \lim_{n\rightarrow\infty} \uparrow  E\left[ (\Kc^{\alpha'_n}_{\tau_2} - \Kc^{\alpha'_n}_{\tau_1})^2\right]^{1/3} 
  \;<\; \infty\;.
  \e*
 Since $\eta'_{\tau_1}$ is in addition non-negative, it is a.s. bounded. Hence, combining \reff{eq : BSDE de la cond minimale - indep alpha futur} and \reff{eqtamp111}, 
 we obtain for $ \alpha' \in \Ab^{\alpha}_{\tau_1}$ 
\begin{align*}
 0&\le   E_{\tau_1}\left[ \Kc^{\alpha'}_{\tau_2}  \right] - \Kc^{\alpha'}_{\tau_1}
 \\&\le
 {C} \;
 ( \Ec^{g}_{\tau_1,\tau_2} [Y^{\alpha'}_{\tau_2}] - \Yc^{\alpha}_{\tau_1})^{1/3}
 \\&= 
 {C} \;
 ( \Ec^{g}_{\tau_1} [\Phi(M^{\alpha'}_{T})] - \Yc^{\alpha}_{\tau_1})^{1/3}\;, \;\;.
 \end{align*}
 Taking the essential infimum in the above inequality  and appealing to \reff{eq : def bar Y} leads to   \reff{eq : BSDE de la cond minimale - cond de minimalite}.\\

\no {\bf Step 3. The uniqueness property for \reff{eq : BSDE de la cond minimale - borne uniforme}-\reff{eq : BSDE de la cond minimale - BSDE}-\reff{eq : BSDE de la cond minimale - cond de minimalite}-\reff{eq : BSDE de la cond minimale - indep alpha futur}}\\
Let us now consider a family $(\tilde Y^{\alpha},\tilde Z^{\alpha}, \tilde K^{\alpha})_{\alpha \in \Ab_{0}}$  of $\Sb_{2}\x\Hb_{2}\x \Kb_{2}$ satisfying \reff{eq : BSDE de la cond minimale - borne uniforme}-\reff{eq : BSDE de la cond minimale - BSDE}-\reff{eq : BSDE de la cond minimale - cond de minimalite}-\reff{eq : BSDE de la cond minimale - indep alpha futur}. Then, \reff{eq : def bar Y}  together with  \reff{eq : BSDE de la cond minimale - BSDE}-\reff{eq : BSDE de la cond minimale - indep alpha futur} applied to $(\tilde Y^{\alpha}, \tilde Z^{\alpha},\tilde K^{\alpha})_{\alpha \in \Ab_{0}}$  imply via a direct comparison argument  that 
\be\label{eqtemp985}
\Yc^{\alpha}_{t}&=& \essinf_{\alpha'\in \Ab_{t}^{\alpha}} \Ec_{t}^{g}[\Phi(M^{\alpha'}_{T})]
\ge \tilde Y^{\alpha}_{t} \;, \qquad {\alpha \in \Ab_{0}}\;,  \quad 0\le t \le T \;.
\ee
 On the other hand, following the exact same line of arguments as the one developed in Step 2 in order to derive \reff{eqtemp1212},  one easily shows that there exists a $\Sb_2$-uniformly bounded family of processes $(\tilde L^\alpha)_{\alpha\in\Ab_{0}}$ such that  
\b*
 \Ec_{t}^{g}[\Phi(M^{\alpha}_{T})] -  \tilde Y^{\alpha}_{t} &=& 
 E_{t}\left[ \int_{t}^{T} \tilde L^{\alpha}_s d\tilde K^{\alpha}_s  \right] 
  \;\le\; 
 C E_{t}\left[ |\tilde K^{\alpha}_{T} - \tilde K^{\alpha}_{t}|^2\right]^{1/2}
 \e*
for all $\alpha\in\Ab_{0},\; 0\le t \le T$,  for some $C>0$. 

Now observe that \reff{eq : BSDE de la cond minimale - cond de minimalite}, applied to $\tilde K^{\alpha}$, and the same arguments as in Lemma \ref{lemma repres Y} provide   the existence of $(\hat \alpha^{n})_{n} \subset \Ab_{t}^{\alpha}$ such that $E_t[\tilde   K^{\hat \alpha^{n}}_{T}-\tilde K^{\alpha}_{t}] \to 0$, $\Pas$  Hence, \reff{eq : BSDE de la cond minimale - borne uniforme} ensures that $E_t[|\tilde   K^{\hat \alpha^{n}}_{T}-\tilde K^{\alpha}_{t}|^2] \to 0$. Since \reff{eq : BSDE de la cond minimale - indep alpha futur} implies $(\tilde   Y^{\hat \alpha^{n}}_t, \tilde   K^{\hat \alpha^{n}}_t)=(\tilde   Y^{  \alpha}_t, \tilde   K^{  \alpha}_t)$ for $n\in\N$, we deduce 
\b*
 \Ec_{t}^{g}[\Phi(M^{\hat \alpha^n}_{T})] -  \tilde Y^{\alpha}_{t}
   &\le&  C \; E_{t}\left[ |\tilde K^{\hat \alpha^n}_{T} - \tilde K^{\alpha}_{t}|^2\right]^{1/2}  \; \rightarrow_{n\rightarrow\infty} \; 0 \;\;.
 \e*
 Combined with \reff{eqtemp985}, this shows that 
 \b*
\tilde Y^{\alpha}_{t} 
&=&  
\essinf_{\alpha'\in \Ab_{t}^{\alpha}}  \Ec^{g}_{t}[\Phi(M^{\alpha'}_{T})]
 \;=\; \Yc^{\alpha}_{t}
 \;, \qquad {\alpha \in \Hb_{2}}\;,  \quad 0\le t \le T \;.
 \e*

The fact that $(\tilde Z^{\alpha},\tilde K^{\alpha})_{\alpha\in\Ab_{0}}=(\Zc^{\alpha},\Kc^{\alpha})_{\alpha\in\Ab_{0}}$ then follows from the uniqueness of the non-linear Doob-Meyer decomposition of \cite[Theorem 3.3]{Pe99}.
\ep

\section{Appendix}

We  report here some standard  results for Lipschitz BSDEs. The first one can be found in, e.g., Theorem 1.5 in \cite{Par98}. The second one is proved for completeness, and by lack of a good reference.

\begin{Proposition}\label{PropStability} (Stability for {Lipschitz} BSDEs)  Let $(Y^1,Z^1)$ and $(Y^2,Z^2)$ in $\Sb_{2}\times\Hb_{2}$ be solutions on $[0,T]$ of Lipschitz BSDEs associated to parameters $(\xi^1,g^1)$ ad $(\xi^2,g^2)$. 
Then the following stability result holds:
\begin{align*}
&\left\| Y^1-Y^2 \right\|_{\Sb_{2}}^2 + \left\| Z^1-Z^2 \right\|_{\Hb_{2}}^2 \\&\le C\left(  \left\| \xi^1-\xi^2 \right\|_{\Lb_2}^2  +  {\int_0^T E\left| g^1-g^2\right|^2(t,Y^1_t,Z^1_t)dt} \right) \;,
\end{align*}
for some constant $C>0$ depending only on $T$ and on the Lipschitz constants of {$g^{1}$ and $g^{2}$}.

\end{Proposition}

 \begin{Proposition} \label{prop: majo bsde appendix} Let the conditions \Hg~hold. Then: 
\begin{enumerate}[(i)]
\item  There exists $C>0$ which only depends on $K_{g}$ and $T$ such that 
 $$
 {\esssup}_{\xi\in \Lb_{0}([0,1])}|\Eg{t}{\xi}|\le  C( 1+ E_{t}\left[|\chi_{g}|^{2}\right]^{\frac12})\;, \qquad 0\le t\le T \;.
 $$

 \item 
 For some $\xi \in \Lb_{2}$ and $t\in [0,T]$, consider a family $(\xi^{\eps})_{\eps\ge 0}\subset \Lb_{0}(\R^{d})$ satisfying $|\xi^{\eps}|\le \xi$ and $\xi^\eps \in L^0(\Fc_{(t+\eps)\wedge T})$, for any $\eps>0$. 
 Then, there exists a family  $(\eta_{\eps})_{\eps>0}\subset \Lb_{0}(\R)$ which converges to $0$ $\Pas$ as $\eps\to 0$ such that 
 $$
{| \Eg{t,t+\eps}{\xi^{\eps}}- E_{t}\left[\xi^{\eps}\right]|}\le \eta_{\eps},\quad \forall \eps\in [0,T-t]. 
 $$
 \item 
 Let $(\xi^\eps)_{\eps>0}$ and $t\in [0,T]$ be as in (ii). 
 Then, there exists a family  $(\eta_{\eps})_{\eps>0}\subset \Lb_{0}(\R)$ which converges to $0$ $\Pas$ as $\eps\to 0$ such that 
 $$
{| \Eg{t-\eps,t}{\xi^{\eps}}- E_{t}\left[\xi^{\eps}\right]|}\le \eta_{\eps},\quad \forall \eps\in [ 0,t]. 
 $$
\end{enumerate}
 \end{Proposition}
 
 \proof \textbf{a.} We first prove (ii) (property (iii) being similar) using the standard linearization argument. Fix $t\in[0,T]$ and set $ Y^{\eps}:=\Eg{\cdot,t+\eps}{\xi^{\eps}}$. Assumption \Hg~implies  that we can find a family of predictable processes $(\rho^{\eps},\gamma^{\eps})$ with values in $[-K_{g},K_{g}]^{d+1}$ such that 
 $$
 L^{\eps}Y^{\eps}+\int_{t}^{\cdot} L^{\eps}_{r} g(r,0,0)dr
 $$ 
 is a martingale on $[t,t+\eps]$,
 with 
 $$
 L^{\eps}_{s}=1+\int_{t}^{s } \rho^{\eps}_{r}L^{\eps}_{r} dr+\int_{t}^{s } \gamma^{\eps}_{r}L^{\eps}_{r} dW_{r},\;t\le s\le t+\eps. 
 $$ 
 In particular,  
 $$
 \Eg{t,t+\eps}{\xi^{\eps}}= L^{\eps}_tY^{\eps}_{t} =E_{t}\left[ L^{\eps}_{t+\eps}\xi^{\eps}+ \int_{t}^{t+\eps} L^{\eps}_{r} g(r,0,0)dr\right]. 
 $$
Condition \Hg~ and the assumption on $(\xi^{\eps})_{\eps>0}$ thus leads to
\b*
| \Eg{t,t+\eps}{\xi^{\eps}}- E_{t}\left[\xi^{\eps}\right]|\le \eta_{\eps}, 
\e*
in which 
$$
\eta_{\eps}:=E_{t}\left[\xi |L^{\eps}_{t+\eps}-L^{\eps}_{t}|+ \chi_{g} \int_{t}^{t+\eps} L^{\eps}_{r}dr \right].
$$
We have:
\begin{align}
\label{eq:tempprop62}
&|\eta_{\eps}|\leq E_t[|\xi|^2]^{1/2} E_t[|L^{\eps}_{t+\eps}-L^{\eps}_{t}|^2]^{1/2} + E_t[|\chi_{g}|^2]^{1/2} E_t\left[\left\vert\int_{t}^{t+\eps} L^{\eps}_{r}dr\right\vert^2\right]^{1/2}\nonumber\\
&\leq E_t[|\xi|^2]^{1/2} E_t[|L^{\eps}_{t+\eps}-L^{\eps}_{t}|^2]^{1/2} + \eps E_t[|\chi_{g}|^2]^{1/2} E_t\left[\sup_{t\leq s\leq t+\eps} |L^{\eps}_{s}|^2\right]^{1/2}.
\end{align}
In addition, 
\begin{align*}
E_t[|L^{\eps}_{t+\eps}-L^{\eps}_{t}|^2]&\leq C E_t\left[\int_t^{t+\eps} |L_r^\eps|^2 dr \right]\\
&\leq \eps C E_t\left[\sup_{t \leq r \leq t+\eps} |L_r^\eps|^2 \right].\\
\end{align*}
Hence,
\begin{align*}
E_t[|L^{\eps}_{t+\eps}-L^{\eps}_{t}|^2]&\leq \eps C \left(1+E_t\left[\sup_{t \leq r \leq t+\eps} |L_r^\eps-L_t^\eps|^2 \right]\right).
\end{align*}
Since $\gamma^{\eps}$ and $\rho^\eps$ are bounded, the quantity $\sup_{t \leq \tau \leq t+\eps} E_\tau[|L^{\eps}_{t+\eps}-L^{\eps}_{\tau}|^2]$ is uniformly bounded. Plugging back this estimate in \eqref{eq:tempprop62} and recalling that $\sup_{t\in [0,T]}E_t[\xi^2]$ is finite $\Pas$ we get that $E_t[|\xi|^2]^{1/2} E_t[|L^{\eps}_{t+\eps}-L^{\eps}_{t}|^2]^{1/2}$ tends to $0$ uniformly in $t$, $\Pas$ as $\eps$ goes to $0$. The second term of \eqref{eq:tempprop62} can be estimated in the same way. 

\textbf{b.} We now prove (i). Pick any $t\in[0,T]$ and $\xi\in \Lb_{0}([0,1])$. The same arguments as above yield
 \b*
\left|\Eg{t}{\xi}\right|&\le& \left|E_{t}\left[L^{\xi}_{T}\xi+ \int_{t}^{T} L^{\xi}_{r} g(r,0,0)dr \right]\right|
\le   E_{t}\left[|L^{\xi}_{T}|+ T|\chi_{g}|\sup_{r\le T}|L^{\xi}_{r}|dr\right]\;,
 \e*
 where $L^{\xi}$ solves
   $$
 L^{\xi}_{s}=1+\int_{t}^{s } \rho^{\xi}_{r}L^{\xi}_{r} dr+\int_{t}^{s } \gamma^{\xi}_{r}L^{\xi}_{r} dW_{r},\;t\le s\le T  \;,
 $$ 
 for some predictable processes $(\rho^{\xi},\gamma^{\xi})$ with values in $[-K_{g},K_{g}]^{d+1}$. Hence, 
 \b*
 \left|\Eg{t}{\xi}\right|\le   E_{t}\left[|L^{\xi}_{T}|+ T|\chi_{g}|\sup_{t\le r\le T}|L^{\xi}_{r}|dr \right].
 \e*
 Since $(\rho^{\xi},\gamma^{\xi})$ are valued in $[-K_{g},K_{g}]^{d+1}$, standard estimates imply that we can find $C>0$, which only depends on $K_{g}$ such that $E_t\left[\sup_{t\le r\le T}|L^{\xi}_{r}|^{2} \right]\le C^{2}$ $\Pas$ The above leads to 
  \b*
 \left|\Eg{t}{\xi}\right|\le  ( C+TC E_{t}\left[ |\chi_{g}|^{2}\right]^{\frac12})\;,
 \e*
 and the arbitrariness of $\xi\in \Lb_{0}([0,1])$ concludes the proof.
    \ep

\end{document}